\documentclass[11pt]{amsart}
\usepackage{comment}
\usepackage{graphics}
\usepackage{tikz}
\usepackage{tikz-cd}
\usepackage{url}
\usepackage{amsmath}
\usepackage{amsfonts}
\usepackage{amssymb}
\usepackage{mathrsfs}
\usepackage{mathtools}
\usepackage{braket}
\mathtoolsset{showonlyrefs=true}

\usepackage[defaultsups]{newtxtext}
\usepackage{amsmath}
\usepackage{amsfonts}
\usepackage{amssymb}
\usepackage{newtxtext}
\usepackage{newtxmath}

\theoremstyle{definition}
\newtheorem{thm}{\textbf{Theorem}}[section]
\newtheorem{cor}[thm]{\textbf{Corollary}}
\newtheorem{defn}[thm]{\textbf{Definition}}
\newtheorem{prop}[thm]{\textbf{Proposition}}
\newtheorem{lem}[thm]{\textbf{Lemma}}
\newtheorem{remk}[thm]{\textbf{Remark}}
\newtheorem{assumption}[thm]{\textbf{Hypothesis}}

\newtheorem{ex}[thm]{\textbf{Example}}

\newcommand{\intZ}{\ensuremath\mathbb{Z}}
\newcommand{\ratQ}{\ensuremath\mathbb{Q}}

\newcommand{\finF}{\ensuremath\mathbb{F}}
\newcommand{\realR}{\ensuremath\mathbb{R}}
\newcommand{\cpxC}{\ensuremath\mathbb{C}}

\newcommand{\projline}[1]{\ensuremath\mathbb{P}^1}

\newcommand{\finFq}[1]{\finF_{#1}}
\newcommand{\projl}[1]{\mathbb{P}^{1}(#1)}

\newcommand{\seqsp}[2]{\mathscr{S}(#1, #2)}% S(f, R)
\newcommand{\valsp}[2]{{V}(#1, #2)}% V(f, R)
\newcommand{\GR}[1]{G_{R}(#1)}

\newcommand{\GQ}[1]{G_{\ratQ}(#1)}
\newcommand{\GFp}[1]{G_{\finF_p}(#1)}

\newcommand{\red}[1]{{\text{red}}_{#1}}
\newcommand{\Gfp}[2]{G(#1, #2)}% G(f, p)
\newcommand{\Kfp}[2]{K(#1, #2)}% K(f, p)
\newcommand{\pvec}[2]{\mathchoice{%
    \begin{pmatrix}
      #1\\
      #2
    \end{pmatrix}
  }{%
    \left(
      \begin{smallmatrix}
        #1\\
        #2
      \end{smallmatrix}
    \right)}{}{}%
}

\newcommand{\bvec}[2]{\mathchoice{%
    \begin{bmatrix}
      #1\\
      #2
    \end{bmatrix}
  }{%
    \left[
      \begin{smallmatrix}
        #1\\
        #2
      \end{smallmatrix}
    \right]}{}{}%
}
\newcommand{\calF}{\ensuremath{\mathcal{F}}}
\newcommand{\calL}{\ensuremath{\mathcal{L}}}
\newcommand{\AIdeclarationtitle}{Use of generative AI}

\usepackage{braket}
\usepackage[hidelinks]{hyperref}

\title{On Lucas zeta functions and Laxton groups}
\author{Iwao Kimura}
\address{Department of Mathematics, University of Toyama,
3190 Gofuku, Toyama 930-8555, Japan}
\email{iwao@sci.u-toyama.ac.jp}

\begin{document}

\begin{abstract}
  Let $R$ be an integral domain, let $P,Q\in R$ with $Q$ a unit, and put
  $f(t)=t^{2}-Pt+Q$ in the polynomial ring $R[t]$. A two-sided sequence
  $w=(w_{n})$ in $R$ with $w_{n+2}=Pw_{n+1}-Qw_{n}$ corresponds to the
  element $w_{1}-w_{0}t$ of the quotient of $R[t]$ by the ideal generated
  by $f(t)$. For two such sequences $v$ and $w$, the product of this ring
  gives a product $v*w$. The sequences whose element is invertible form,
  modulo scalars, the Laxton group. For real $P$ and $Q$ with
  $P^{2}-4Q>0$, we study how its group structure is reflected in the Lucas
  zeta functions $\zeta_{w}(s)=\sum_{n\ge1}w_{n}^{-s}$ of a complex variable
  $s$, for sequences with $w_{n}>0$ ($n\ge1$). Let $\theta_{1}$ and
  $\theta_{2}$ be the roots of $f$, with $\theta_{1}$ larger than the
  absolute value of $\theta_{2}$, and put $A_{w}=w_{1}-w_{0}\theta_{2}$
  and $B_{w}=w_{1}-w_{0}\theta_{1}$. The ratio $\lambda_{w}=B_{w}/A_{w}$
  gives an isomorphism of the Laxton group over the real numbers onto the
  multiplicative group of nonzero real numbers, and $\zeta_{w}$ depends on
  $w$ only through $A_{w}$ and $\lambda_{w}$. At a negative integer $-m$,
  apart from explicit exceptions, $\zeta_{w}(-m)$ is the value of a binary form of
  degree $m$ at $(A_{w},B_{w})$ plus a correction term. The correction
  term vanishes when $f$ is irreducible over the rationals, but in general
  it makes $\zeta_{w}(-m)$ irrational for rational $P$, $Q$, $w_{0}$,
  $w_{1}$; this corrects a closed formula in the literature. For an odd
  prime $p$, the values attached to the powers $w^{*k}$ satisfy
  congruences modulo powers of $p$ in the exponent $k$, both for primes
  that divide $P^{2}-4Q$ and for those that do not. These values are
  interpolated $p$-adically in $k$. When $Q=\pm1$ and $P$ is an integer, the sequences are the traces of $\alpha\theta_{1}^{n}$
  for elements $\alpha$ of the real quadratic field generated by
  $\theta_{1}$, and $\theta_{1}$ is a power of its fundamental unit.
\end{abstract}

\maketitle{}
\bibliographystyle{amsalpha}

\section{Introduction}
\label{sec:introduction}

For real numbers $P$ and $Q$, the Lucas sequence $(\calF_{n})$ is given by
$\calF_{0}=0$, $\calF_{1}=1$ and $\calF_{n+2}=P\calF_{n+1}-Q\calF_{n}$; for
$(P,Q)=(1,-1)$, $\calF_{n}$ is the $n$-th Fibonacci number. Let $s$ be a
complex variable, and put $x^{-s}=\exp(-s\log x)$ for real $x>0$, with the
real logarithm. The Fibonacci zeta function $\sum_{n\ge1}\calF_{n}^{-s}$
was continued meromorphically to the whole complex plane by
Egami~\cite{MR1744904} and Navas~\cite{MR1866357}. Kamano~\cite{MR3064567}
did the same for the Lucas sequences, and for their twists by Dirichlet
characters. An early case with general
Binet coefficients is the series $\sum_{n\ge1}(ar^{n}+br^{-n})^{-s}$,
with real numbers $r>1$, $a$ and $b$~\cite{SilvermanBradleyDarling1999}.
Smajlovi\'c, \v{S}abanac and \v{S}\'ceta studied Hurwitz-type and
Tribonacci analogs~\cite{zbMATH07676374,zbMATH07676373}. Meher and Rout
continued multiple Lucas zeta functions and their twists by additive and
Dirichlet characters~\cite{MR3852570,MR4333647}. For integer linear
recurrences of arbitrary degree whose minimal polynomial is irreducible
and has a dominant root, Serrano Holgado and Navas~\cite{MR4594520}
determined the poles and residues and proved that the values at negative
integers are rational; the Lerch-type series are treated
in~\cite{zbMATH08107325}.

Let $\varphi=(1+\sqrt5)/2$ be the golden ratio, with $\sqrt5>0$, and let
$\mathrm{Tr}$ denote the trace from a real quadratic field to the field of
rational numbers. The Lucas numbers are the traces
$\mathrm{Tr}(\varphi^{n})$. More
generally, one can study the traces $\mathrm{Tr}(\alpha\varepsilon^{n})$,
where $\varepsilon$ is the fundamental unit of a real quadratic field and
$\alpha$ is an element of that field. Assaf, Kuan, Lowry-Duda and
Walker~\cite{AKLW2024} introduced the sequences
$\mathrm{Tr}(\varepsilon^{n})$ and $\mathrm{Tr}(\varepsilon^{n}/\sqrt{D})$,
where $D>0$ is the discriminant of the field, and studied the zeta function of
the latter. Sequences of this type attached to cubic fields that are not
totally real are studied in~\cite{Kimura2025MSJ}.

In this paper $\alpha$ is arbitrary. We write $\intZ$, $\ratQ$, $\realR$
and $\cpxC$ for the integers, the rational, the real and the complex
numbers, and $X^{\times}$ for the group of units of a ring $X$. Let $R$
be an integral domain, let $P, Q\in R$ with $Q\in R^{\times}$, and put
$f(t)=t^{2}-Pt+Q$ in the polynomial ring $R[t]$. A sequence
$w=(w_{n})_{n\in\intZ}$ in $R$ with $w_{n+2}=Pw_{n+1}-Qw_{n}$ corresponds to
the element $w_{1}-w_{0}t$ of the quotient ring $R[t]/(f(t))$ of $R[t]$ by
the ideal generated by $f(t)$. For two such sequences $v$ and $w$, the
product of this ring gives a product $v*w$. The sequences whose element is
invertible form, modulo multiplication by $R^{\times}$, the Laxton group
$\GR{f}$~\cite{zbMATH03358544,zbMATH03358545,MR3902117}. We write
$\bvec{w_{1}}{w_{0}}$ for the class of $w$, and $w^{*k}$ for the $k$-th
power of $w$ under $*$ ($k\in\intZ$). We study how the group structure of
$\GR{f}$ is reflected in the Lucas zeta functions
$\zeta_{w}(s)=\sum_{n\ge1}w_{n}^{-s}$.

Let $P, Q\in\realR$, put $d=P^{2}-4Q$, and assume $d>0$. Let
$\theta_{1}$ and $\theta_{2}$ be the roots of $f$, with $\theta_{1}$
larger than the absolute value $\lvert\theta_{2}\rvert$ of $\theta_{2}$.
Put $A_{w}=w_{1}-w_{0}\theta_{2}$ and $B_{w}=w_{1}-w_{0}\theta_{1}$; up
to the factors $1/\sqrt d$ and $-1/\sqrt d$, these are the coefficients of
Binet's formula. When $P\in\intZ$ and $Q=\pm1$, the
integer $d$ is not a square, and $\theta_{1}=\varepsilon^{g}$ for the fundamental
unit $\varepsilon$ of the real quadratic field $\ratQ(\sqrt d)$ and some
integer $g\ge1$, and the sequences $w$ are
exactly the sequences $\mathrm{Tr}(\alpha\varepsilon^{gn})$
(Proposition~\ref{prop:units}). We assume Kamano's condition $(C)$ on
$(P,Q)$, stated in Definition~\ref{def:condc}, and that the sequences
involved are positive, as stated in Hypothesis~\ref{hyp:positivity}. The
main results are as follows.
\begin{enumerate}
\item The map $\bvec{w_{1}}{w_{0}}\mapsto\lambda_{w}=B_{w}/A_{w}$ is an
  isomorphism $G_{\realR}(f)\to\realR^{\times}$; this is the case
  $R=\realR$ of an isomorphism of Aoki and Kida
  (Proposition~\ref{prop:lambda}). Up to the factor $d^{s/2}A_{w}^{-s}$
  and finitely many terms, $\zeta_{w}(s)$ is the value at $\lambda_{w}$ of
  a function $\Phi_{n_{1}}(s;\lambda)$ that depends only on $f$ and on a
  positive integer $n_{1}$; it is defined in~\eqref{eq:Phi}
  (Theorem~\ref{thm:expforseq}). Hence $\zeta_{w}$ depends on the class
  of $w$ only through $A_{w}$ and $\lambda_{w}$, and $\zeta_{v*w}$ is
  obtained by evaluating the same function at $\lambda_{v}\lambda_{w}$
  (Corollary~\ref{cor:zetaprod}).
\item Let $m\ge1$, and assume that $\theta_{1}^{m-i}\theta_{2}^{i}\neq1$ for
  $0\le i\le m$; then $s=-m$ is not a pole. Moreover
  $\zeta_{w}(-m)=H_{m}(A_{w},B_{w})+E_{w}(m)$, where $H_{m}$ is a binary
  form of degree $m$ that depends only on $f$ and $m$, defined in
  Corollary~\ref{cor:negform}, and $E_{w}(m)$ is the term given in
  Theorem~\ref{thm:negval}. The term $E_{w}(m)$ vanishes when $f$ is
  irreducible over $\ratQ$, but not in general. Then $\zeta_{w}(-m)$ can
  be irrational for rational $P$, $Q$, $w_{0}$, $w_{1}$
  (Remark~\ref{remk:negval}). This corrects a closed formula of
  Kamano~\cite[Theorem~7]{MR3064567}. It also shows that the rationality
  result of~\cite{MR4594520} does not extend to separable reducible
  minimal polynomials, as suggested in~\cite[\S7]{MR4594520}. When
  $E_{w}(m)=0$ and $P$, $Q$, $w_{0}$, $w_{1}$ are rational,
  $\zeta_{w}(-m)$ is rational (Theorem~\ref{thm:negrational}), and
  $\zeta_{v*w}(-m)=H_{m}(A_{v}A_{w},B_{v}B_{w})$.
\item Let $P, Q\in\intZ$, let $d$ be not a square, and put
  $K=\ratQ(\sqrt d)$. Let $p$ be an odd prime, and let
  $\delta_{p}\in\set{1,-1,0}$ be the Legendre symbol $(d/p)$. Put
  $\tau_{w}=1$, $\mathrm{N}(A_{w})$ or $w_{1}-w_{0}P/2$ according as
  $\delta_{p}=1$, $-1$ or $0$; here $\mathrm{N}$ is the norm from $K$ to
  $\ratQ$. Let $w_{0}, w_{1}\in\intZ$ and $m\ge1$. Assume that
  $p\nmid\mathrm{N}(A_{w})$, that
  $p\nmid\mathrm{N}(1-\theta_{1}^{m-i}\theta_{2}^{i})$ for $0\le i\le m$,
  and that $p\nmid Q$ if $\delta_{p}=0$. Then
  \[
    \zeta_{w^{*(p-\delta_{p})}}(-m)\equiv
    \tau_{w}^{m}\,\zeta_{\calF}(-m)\pmod{p}
  \]
  (Corollary~\ref{cor:fermat}). This is a congruence between rational
  numbers whose denominators are prime to $p$. Here $\calF$ is the Lucas
  sequence $(\calF_{n})$, whose class is the identity of $\GQ{f}$. For
  $p\mid d$ the coefficients of $H_{m}$ contain $d^{-m/2}$, but the
  values are still $p$-integral. More generally, let $e\ge0$, and let
  $k, k'\in\intZ$, possibly negative, with
  $k\equiv k'\pmod{(p-\delta_{p})p^{e}}$. Then
  $H_{m}(A_{w}^{k},B_{w}^{k})$ and $H_{m}(A_{w}^{k'},B_{w}^{k'})$ are
  congruent modulo $p^{e+1}$, up to a power of $\tau_{w}$
  (Theorem~\ref{thm:congruence}). These values are
  $\zeta_{w^{*k}}(-m)$ and $\zeta_{w^{*k'}}(-m)$ when $w^{*k}$ and
  $w^{*k'}$ satisfy Hypothesis~\ref{hyp:positivity}
  (Remark~\ref{remk:interppos}). Put $M=p-1$,
  $p^{2}-1$ or $p(p-1)$ according as $\delta_{p}=1$, $-1$ or $0$, and
  put $k=a+Mn$ with $a$ fixed. Then $H_{m}(A_{w}^{k},B_{w}^{k})$ is a
  $p$-adic Lipschitz function of $n\in\intZ$
  (Corollaries~\ref{cor:period} and~\ref{cor:padicinterp}). These
  congruences concern the exponent $k$ in the Laxton group. They follow
  from the form $H_{m}$ of item~2 and from the congruence
  $w^{*(p-\delta_{p})p^{e}}\equiv\tau_{w}^{p^{e}}\calF\pmod{p^{e+1}}$ of
  sequences (Lemma~\ref{lem:laxtonpower}). If $p\nmid Q$, the exponent
  $p-\delta_{p}$ is the order of the group $\Kfp{f}{p}/\Gfp{f}{p}$ of
  Aoki--Kida~\cite{MR3902117}. Congruences in $m$ of Kummer
  type are not treated here; we will treat them in a separate paper.
\end{enumerate}

For integer sequences with irreducible $f$, the meromorphic continuation
and the rationality of the values at negative integers are special cases
of~\cite{MR4594520}, and the continuation of the partial zeta functions
and $L$-functions in \S\ref{sec:lucas-l-functions} follows
from~\cite{zbMATH08107325}; for $w=\calF$ the latter are due to
Kamano~\cite[\S4]{MR3064567}. Serrano Holgado and Navas remark that their
proofs extend to real coefficients~\cite[\S7]{MR4594520}. We give
complete proofs, because we need the explicit form in terms of $A_{w}$
and $B_{w}$, and because the reducible case needs the term $E_{w}(m)$.

In Section~\ref{sec:laxton-groups} we recall the Laxton group. In
Section~\ref{sec:lucas-zeta-functions} we treat the continuation, the
poles and residues, products, partial zeta functions and $L$-functions.
In Section~\ref{sec:special-values-lucas} we treat the values at negative
integers and the congruences. In Section~\ref{sec:units} we treat the
case of units of real quadratic fields.

\subsection*{Notation}
We write $\intZ_{\ge0}$ and $\intZ_{\ge1}$ for the nonnegative and the
positive integers. For $m\in\intZ_{\ge0}$, $m!$ is the factorial.
For integers $a$, $b$ and a prime $p$,
$a\mid b$ means that $a$ divides $b$, $a\nmid b$ that it does not, and
$p^{e}\parallel b$ that $p^{e}\mid b$ and $p^{e+1}\nmid b$.

For a ring $X$, $X[t]$ is the polynomial ring,
$X^{2} = X\times X$ the product, $X^{\intZ}$ the set of sequences $(x_{n})_{n\in\intZ}$ in $X$,
and $\set{x | \dots}$ a set.
For an element $x$ of a commutative ring
$X$, $(x)=xX$ is the principal ideal; for an ideal $I$ of $X$,
$y\equiv z\pmod{I}$ means $y-z\in I$, and $\pmod{x}$ means $\pmod{xX}$.

For $z\in\cpxC$, $\lvert z\rvert$ is the absolute value,
$\Re z$ the real part, $\exp z$ the exponential,
and $\sqrt{-1}$ a fixed square root of $-1$; for real $x\ge0$, $\sqrt x\ge0$;
$\log$ is the real logarithm of a positive real number, and
$\pi=3.14\ldots$. For a function $\psi$ of $s$, $\frac{\mathrm d}{\mathrm ds}\psi$
is its derivative. $\mathop{\mathrm{Res}}_{s=s_{0}}$ denotes the residue at $s=s_{0}$, and
$O(\cdot)$ is the Landau symbol.

For a prime $p$, $\finFq{p}$ is
the field $\intZ/p\intZ$ and $\mathbb F_{p^{2}}$ the field with $p^{2}$ elements, $v_{p}$ the $p$-adic valuation on $\ratQ$,
$\lvert x\rvert_{p}=p^{-v_{p}(x)}$ the $p$-adic absolute value,
$\intZ_{p}$ the ring of $p$-adic integers, to which $v_{p}$ and
$\lvert\cdot\rvert_{p}$ extend, and
$\intZ_{(p)}=\set{x\in\ratQ | v_{p}(x)\ge0}$; for odd $p$ and an integer
$d$, $(d/p)\in\set{1,-1,0}$ is the Legendre symbol, with $(d/p)=0$ if
$p\mid d$.

For $d\in\intZ$ not a square,
$\ratQ(\sqrt d)$ is the quadratic field generated by a square root of
$d$. For a quadratic field $K$, $\mathcal O$ is its ring of integers,
$d_{K}$ its discriminant, $\sigma$ its nontrivial automorphism,
$\mathrm{Tr}_{K/\ratQ}(x)=x+\sigma(x)$ and
$\mathrm{N}(x)=\mathrm{N}_{K/\ratQ}(x)=x\sigma(x)$ for $x\in K$. For
$\theta\in\mathcal O$, $\intZ[\theta]$ is the subring generated by
$\theta$, and $[\mathcal O:\intZ[\theta]]$ the index. A prime ideal
$\mathfrak p$ of $\mathcal O$ is \emph{above} a prime $p$ if
$\mathfrak p\cap\intZ=p\intZ$; then
$\mathcal O_{\mathfrak p}=\set{x/y | x,y\in\mathcal O,\ y\notin\mathfrak p}$
is the localization. It is a subring of $K$ containing $\mathcal O$, and
every $y\in\mathcal O\setminus\mathfrak p$ is a unit of
$\mathcal O_{\mathfrak p}$. In particular, every $x\in\mathcal O$ with
$p\nmid\mathrm{N}(x)$ is a unit of $\mathcal O_{\mathfrak p}$, since
$x\in\mathfrak p$ would give
$\mathrm{N}(x)=x\sigma(x)\in\mathfrak p\cap\intZ=p\intZ$. We have
$\mathcal O_{\mathfrak p}\cap\ratQ=\intZ_{(p)}$: integers prime to $p$
are not in $\mathfrak p$; conversely, if
$y\in\mathcal O_{\mathfrak p}\cap\ratQ$ had $v_{p}(y)<0$, then $1/p$
would lie in $\mathcal O_{\mathfrak p}$, say $1/p=x/y'$ with
$y'\notin\mathfrak p$, but then $y'=px\in\mathfrak p$. Hence, for
$y\in\ratQ$ and $h\ge0$, $y\in p^{h}\mathcal O_{\mathfrak p}$ implies
$y\in p^{h}\intZ_{(p)}$. For a field
extension, $\mathrm{Gal}$ denotes its Galois group, and $\mathrm{id}$ is
the identity map.

\section{Laxton groups}
\label{sec:laxton-groups}

We follow Aoki--Kida~\cite{MR3902117}, who give a detailed account of
Laxton groups.

Let $R$ be an integral domain, and let $P, Q\in R$ with $Q\in R^{\times}$.
Put $f(t)=t^{2}-Pt+Q$ and
\begin{equation}
  \label{eq:seqsp}
  \seqsp{f}{R}\coloneqq\set{(w_{n})_{n\in\intZ}\in R^{\intZ} | w_{n+2}=Pw_{n+1}-Qw_{n}\ \text{for all $n$}}.
\end{equation}
Since $Q\in R^{\times}$, the recurrence can be run backwards,
$w_{n}=Q^{-1}(Pw_{n+1}-w_{n+2})$. Hence a sequence in $\seqsp{f}{R}$ is
determined by $(w_{1},w_{0})\in R^{2}$, and every pair occurs. Under
termwise operations, the maps
\begin{equation}
  \label{eq:three-modules}
  \begin{aligned}
    \seqsp{f}{R}&\to\valsp{f}{R}\coloneqq R^{2}\to R[t]/(f(t)),\\
    (w_{n})&\mapsto\pvec{w_{1}}{w_{0}}\mapsto w_{1}-w_{0}t+(f(t)),
  \end{aligned}
\end{equation}
are isomorphisms of $R$-modules. We transport the ring structure of
$R[t]/(f(t))$ to $\seqsp{f}{R}$ and $\valsp{f}{R}$, and write $*$ for the
product. The identity element is the sequence with $w_{0}=0$, $w_{1}=1$.

\begin{ex}\label{ex:sequences}
  \begin{enumerate}
  \item For $f(t)=t^{2}-t-1$, the identity of $\seqsp{f}{\intZ}$ is the
    Fibonacci sequence $(\calF_{n})$, with $\calF_{-1}=1$, $\calF_{0}=0$,
    $\calF_{1}=\calF_{2}=1$.
  \item For $f(t)=t^{2}-3t+2=(t-1)(t-2)$, the sequence $\mathcal{M}_{n}=2^{n}-1$ lies
    in $\seqsp{f}{\ratQ}$, since $3\mathcal{M}_{n+1}-2\mathcal{M}_{n}=2^{n+2}-1$. It is the
    identity ($\mathcal{M}_{0}=0$, $\mathcal{M}_{1}=1$; note $\mathcal{M}_{-1}=-1/2$).
  \end{enumerate}
\end{ex}

\begin{defn}\label{def:laxton}
  The \emph{Laxton group} is $\GR{f}\coloneqq\valsp{f}{R}^{\times}/R^{\times}$.
  The class of $\pvec{w_{1}}{w_{0}}$ is written $\bvec{w_{1}}{w_{0}}$, and
  the product in $\GR{f}$ is again written $*$. For $w\in\seqsp{f}{R}$ and
  $k\ge0$, $w^{*k}$ is the $k$-th power of $w$ under $*$; for $k<0$ it is
  defined when $w$ is invertible.
\end{defn}

\begin{lem}\label{lem:multiplication}
  For $v, w\in\seqsp{f}{R}$,
  \[
    v*w=\pvec{w_{1}v_{1}-Qw_{0}v_{0}}{w_{1}v_{0}+w_{0}v_{1}-Pw_{0}v_{0}}.
  \]
\end{lem}
\begin{proof}
  In $R[t]/(f(t))$ we have $t^{2}=Pt-Q$, so
  \[
    (v_{1}-v_{0}t)(w_{1}-w_{0}t)
    =(v_{1}w_{1}-Qv_{0}w_{0})-(v_{1}w_{0}+v_{0}w_{1}-Pv_{0}w_{0})\,t.
    \qedhere
  \]
\end{proof}

\begin{defn}\label{def:lucasseq}
  The \emph{Lucas sequence} $\calF\in\seqsp{f}{R}$ has $\calF_{0}=0$,
  $\calF_{1}=1$; it is the identity element. The \emph{companion Lucas
  sequence} $\calL\in\seqsp{f}{R}$ has $\calL_{0}=2$, $\calL_{1}=P$.
\end{defn}

\begin{prop}\label{prop:binet}
  Let $F$ be a field of characteristic $\neq2$, let $P, Q\in F$ with
  $Q\neq0$ and $d\coloneqq P^{2}-4Q\neq0$. Fix a square root $\sqrt d$ in an
  algebraic closure $F^{\mathrm{alg}}$, and put $\theta_{1,2}=(P\pm\sqrt d)/2$. For
  $w\in\seqsp{f}{F}$ put $A=A_{w}=w_{1}-w_{0}\theta_{2}$ and
  $B=B_{w}=w_{1}-w_{0}\theta_{1}$. Then
  \begin{equation}
    \label{eq:binet}
    w_{n}=\frac{A\theta_{1}^{n}-B\theta_{2}^{n}}{\sqrt d}\qquad(n\in\intZ).
  \end{equation}
\end{prop}
\begin{proof}
  We have $\theta_{1}+\theta_{2}=P$, $\theta_{1}\theta_{2}=Q\neq0$ and
  $\theta_{1}-\theta_{2}=\sqrt d$. Since $\theta_{i}^{2}=P\theta_{i}-Q$,
  the right-hand side of \eqref{eq:binet} satisfies the recurrence, so it
  suffices to check $n=0,1$. At $n=0$ it is $(A-B)/\sqrt d=w_{0}$. At
  $n=1$ it is $(w_{1}\theta_{1}-w_{0}Q-w_{1}\theta_{2}+w_{0}Q)/\sqrt d=w_{1}$.
\end{proof}

\begin{prop}\label{prop:prodbinet}
  In the setting of Proposition~\ref{prop:binet}, let $v, w\in\seqsp{f}{F}$.
  Then $A_{v*w}=A_{v}A_{w}$ and $B_{v*w}=B_{v}B_{w}$. Hence
  \[
    (v*w)_{n}=\frac{A_{v}A_{w}\theta_{1}^{n}-B_{v}B_{w}\theta_{2}^{n}}{\sqrt d}
    \quad\text{and}\quad
    (w^{*k})_{n}=\frac{A_{w}^{k}\theta_{1}^{n}-B_{w}^{k}\theta_{2}^{n}}{\sqrt d}
  \]
  for $k\ge0$, and for all $k\in\intZ$ if $A_{w}B_{w}\neq0$.
\end{prop}
\begin{proof}
  The map $F[t]/(f(t))\to F^{\mathrm{alg}}$, $t\mapsto\theta_{2}$, is a ring
  homomorphism sending $w_{1}-w_{0}t$ to $A_{w}$; likewise $t\mapsto\theta_{1}$
  sends it to $B_{w}$. By the definition of $*$, this gives
  $A_{v*w}=A_{v}A_{w}$ and $B_{v*w}=B_{v}B_{w}$. The formulas follow from
  Proposition~\ref{prop:binet}. For negative $k$ note that
  $w\in\seqsp{f}{F}^{\times}$ when $A_{w}B_{w}\neq0$, since
  $F[t]/(f)$ is either a field or $F\times F$ (as $d\neq0$).
\end{proof}

\begin{ex}\label{ex:fibonacci}
  For the companion Lucas sequence, $A_{\calL}=P-2\theta_{2}=\sqrt d$ and
  $B_{\calL}=P-2\theta_{1}=-\sqrt d$. By Proposition~\ref{prop:prodbinet},
  $\calL^{*2}=d\,\calF$. Hence, if $d\in R^{\times}$, the class of $\calL$
  has order dividing $2$ in $\GR{f}$. For example, for $f(t)=t^{2}-t-1$
  and $R=\ratQ$ we have $d=5$ and $\bvec{1}{2}^{*2}=\bvec{5}{0}=\bvec{1}{0}$.
  Over $R=\intZ$, the element $1-2t$ has norm $A_{\calL}B_{\calL}=-5$,
  which is not a unit, so $\calL$ has no class in $G_{\intZ}(f)$.
\end{ex}

\section{Lucas zeta functions}
\label{sec:lucas-zeta-functions}

\subsection{Definition and the condition \texorpdfstring{$(C)$}{(C)}}
\label{sec:lucas-zeta-l}

In this section $R=\realR$ and $f(t)=t^{2}-Pt+Q\in\realR[t]$ with $Q\neq0$.
When $d>0$, $\sqrt d$ denotes the positive square root. We use the
notation of Section~\ref{sec:laxton-groups}: the module $\seqsp{f}{R}$
of~\eqref{eq:seqsp}, the product $*$ and the powers $w^{*k}$
(Definition~\ref{def:laxton}), the Lucas sequence $\calF$
(Definition~\ref{def:lucasseq}), and $\theta_{1}$, $\theta_{2}$, $A_{w}$,
$B_{w}$ (Proposition~\ref{prop:binet}).

\begin{defn}\label{def:zeta}
  For $w\in\seqsp{f}{\realR}$ the \emph{Lucas zeta function} is
  \[
    \zeta_{w}(s)\coloneqq\sum_{n=1}^{\infty}w_{n}^{-s}
  \]
  for those $s\in\cpxC$ for which the series converges; here $w_{n}>0$ is assumed and
  $w_{n}^{-s}=\exp(-s\log w_{n})$ with $\log w_{n}\in\realR$. Conditions for
  this are given in Lemma~\ref{lem:convergence} and
  Hypothesis~\ref{hyp:positivity} below.
\end{defn}

\begin{remk}\label{remk:scaling}
  If $w'=aw$ with $a>0$, then $\zeta_{w'}(s)=a^{-s}\zeta_{w}(s)$. Thus,
  up to this factor, $\zeta_{w}$ depends only on the class of $w$ modulo
  positive scalars. For $a<0$, $\zeta_{w'}$ is not defined.
\end{remk}

\begin{defn}\label{def:condc}
  We say that $P, Q\in\realR$ satisfy the condition $(C)$ if
\begin{equation}
  \label{eq:hypothesis}
  P > 0,\quad Q\neq 0,\quad\text{and}\quad
  \begin{dcases}
    Q<P-1 & \text{if $0<P\le2$},\\
    Q\le P-1 & \text{if $2<P$}.
  \end{dcases}
\end{equation}
\end{defn}

\begin{lem}[{\cite[Proposition~3]{MR3064567}}]\label{lem:kamano}
  Under the condition $(C)$, the following two assertions hold:
  \begin{enumerate}
  \item $d = P^2-4Q>0$, $\theta_{1}>1$,
    $\theta_{1} > \lvert \theta_{2}\rvert>0$,
  \item
    $\calF_{n}=
    \dfrac{\theta_{1}^{n}-\theta_{2}^{n}}{\theta_{1}-\theta_{2}} > 0$,
    for all $n\in\intZ_{\ge1}$.
  \end{enumerate}
\end{lem}

\begin{lem}\label{lem:shift}
  Let $F$ be a field, $f(t)=t^{2}-Pt+Q\in F[t]$, $Q\neq0$, and
  $v, w\in\seqsp{f}{F}$. Put $u=v*w$, i.e.\
  $u_{1}=w_{1}v_{1}-Qw_{0}v_{0}$ and
  $u_{0}=w_{1}v_{0}+w_{0}v_{1}-Pw_{0}v_{0}$
  (Lemma~\ref{lem:multiplication}). Then
  \begin{equation}
    \label{eq:shift}
    u_{n} = w_{1}v_{n} - Qw_{0}v_{n-1}\qquad(n\in\intZ).
  \end{equation}
  In particular, taking $v=\calF$ (the identity element),
  $w_{n}=w_{1}\calF_{n}-Qw_{0}\calF_{n-1}$.
\end{lem}
\begin{proof}
  Both sides of \eqref{eq:shift} satisfy the recurrence
  $x_{n+2}=Px_{n+1}-Qx_{n}$ in $n$, so it suffices to compare them at
  $n=0,1$. At $n=1$ the right-hand side is $w_{1}v_{1}-Qw_{0}v_{0}=u_{1}$.
  At $n=0$, using $v_{-1}=(Pv_{0}-v_{1})/Q$, it is
  $w_{1}v_{0}-w_{0}(Pv_{0}-v_{1})=u_{0}$. The last assertion is
  \eqref{eq:shift} with $v=\calF$, since $\calF*w=w$.
\end{proof}

\begin{prop}\label{prop:trace}
  Let $P, Q\in\ratQ$, $Q\neq0$, and assume that $d=P^{2}-4Q$ is not a
  square in $\ratQ$. Put $K=\ratQ(\sqrt{d})$ and let $\sigma$ be the
  nontrivial automorphism of $K$. For $w\in\seqsp{f}{\ratQ}$ put
  $A_{w}=w_{1}-w_{0}\theta_{2}\in K$. Then
  \begin{enumerate}
  \item $\sigma(A_{w})=B_{w}$, and
    $w_{n}=\mathrm{Tr}_{K/\ratQ}\bigl(A_{w}\theta_{1}^{n}/\sqrt{d}\bigr)$
    for all $n\in\intZ$;
  \item $A_{v*w}=A_{v}A_{w}$ for $v, w\in\seqsp{f}{\ratQ}$; hence
    $(v*w)_{n}=\mathrm{Tr}_{K/\ratQ}\bigl(A_{v}A_{w}\theta_{1}^{n}/\sqrt{d}\bigr)$.
  \end{enumerate}
\end{prop}
\begin{proof}
  (1) By Proposition~\ref{prop:binet},
  $w_{n}=(A\theta_{1}^{n}-B\theta_{2}^{n})/\sqrt{d}$ with $A=A_{w}$ and
  $B=w_{1}-w_{0}\theta_{1}$. As $w_{0}, w_{1}\in\ratQ$ and
  $\sigma(\theta_{2})=\theta_{1}$, we have $\sigma(A)=B$; also
  $\sigma(\theta_{1})=\theta_{2}$ and $\sigma(\sqrt{d})=-\sqrt{d}$. Hence
  $\sigma(A\theta_{1}^{n}/\sqrt{d})=-B\theta_{2}^{n}/\sqrt{d}$, and the sum
  of the two conjugates is $w_{n}$.
  (2) The first assertion is Proposition~\ref{prop:prodbinet}; then apply
  (1) to $v*w$.
\end{proof}

\begin{lem}\label{lem:positivity}
  Assume the condition $(C)$, and let $w\in\seqsp{f}{\realR}$ satisfy
  \[
    w_{0}\ge0,\qquad w_{1}>0,\qquad A=w_{1}-w_{0}\theta_{2}\ge0.
  \]
  Then $w_{n}>0$ for all $n\ge1$. If $Q<0$, the first two conditions
  imply $A>0$.
\end{lem}
\begin{proof}
  By Lemma~\ref{lem:kamano}, $\calF_{n}>0$ for $n\ge1$, $\calF_{0}=0$,
  and $\theta_{1}>1$. By Lemma~\ref{lem:shift},
  $w_{n}=w_{1}\calF_{n}-Qw_{0}\calF_{n-1}$.

  If $Q<0$, then $w_{n}\ge w_{1}\calF_{n}>0$ for $n\ge1$. Also
  $\theta_{2}=Q/\theta_{1}<0$, so $A=w_{1}+w_{0}\lvert\theta_{2}\rvert>0$.

  If $Q>0$, then $\theta_{2}=Q/\theta_{1}>0$. We have $w_{1}>0$. For
  $n\ge2$, the identity
  $\calF_{n}=\theta_{1}\calF_{n-1}+\theta_{2}^{n-1}$ (multiply by
  $\theta_{1}-\theta_{2}\neq0$ to check it) gives
  $\calF_{n}>\theta_{1}\calF_{n-1}$, hence
  $Q\calF_{n-1}=\theta_{1}\theta_{2}\calF_{n-1}<\theta_{2}\calF_{n}$. If
  $w_{0}=0$ then $w_{n}=w_{1}\calF_{n}>0$; if $w_{0}>0$ then
  $w_{n}>w_{1}\calF_{n}-w_{0}\theta_{2}\calF_{n}=A\calF_{n}\ge0$.
\end{proof}

\begin{lem}\label{lem:convergence}
  Assume the condition $(C)$, and let $w\in\seqsp{f}{\realR}$ satisfy
  $w_{n}>0$ for all $n\ge1$ and $A=w_{1}-w_{0}\theta_{2}>0$. Then
  $\zeta_{w}(s)=\sum_{n\ge1}w_{n}^{-s}$ (with $w_{n}^{-s}=\exp(-s\log w_{n})$,
  $\log w_{n}\in\realR$) converges absolutely and uniformly on compact
  subsets of $\Re s>0$, and defines a holomorphic function there.
\end{lem}
\begin{proof}
  By Lemma~\ref{lem:kamano}, $\theta_{1}>1$ and
  $r\coloneqq\lvert\theta_{2}\rvert/\theta_{1}<1$. By
  Proposition~\ref{prop:binet},
  $w_{n}=\theta_{1}^{n}\bigl(A-B(\theta_{2}/\theta_{1})^{n}\bigr)/\sqrt{d}$,
  and $\lvert B(\theta_{2}/\theta_{1})^{n}\rvert=\lvert B\rvert r^{n}$. Choose $n_{1}$ with $\lvert B\rvert r^{n}\le A/2$
  for $n\ge n_{1}$; then $w_{n}\ge c\,\theta_{1}^{n}$ for $n\ge n_{1}$, where
  $c=A/(2\sqrt{d})>0$. Let $\Omega$ be a compact subset of $\Re s>0$, and put
  $\tau_{0}=\min_{\Omega}\Re s>0$ and $\tau_{1}=\max_{\Omega}\Re s$. For $s\in\Omega$
  and $n\ge n_{1}$,
  $\lvert w_{n}^{-s}\rvert=w_{n}^{-\Re s}
  \le\max(c^{-\tau_{0}},c^{-\tau_{1}})\,\theta_{1}^{-n\tau_{0}}$,
  a convergent geometric series. The Weierstrass $M$-test gives the
  assertion; the finitely many terms with $n<n_{1}$ are entire in $s$
  because $w_{n}>0$.
\end{proof}

\begin{prop}\label{prop:prodpositive}
  Assume the condition $(C)$ with $Q<0$. Let $v, w\in\seqsp{f}{\realR}$
  satisfy $v_{0}, w_{0}\ge0$ and $v_{1}, w_{1}>0$. Then $(v*w)_{n}>0$ for
  all $n\ge1$, and $\zeta_{v*w}(s)$ converges as in
  Lemma~\ref{lem:convergence}.
\end{prop}
\begin{proof}
  By Lemma~\ref{lem:positivity}, $v_{n}>0$ for $n\ge1$, and $v_{0}\ge0$.
  By Lemma~\ref{lem:shift},
  $(v*w)_{n}=w_{1}v_{n}+\lvert Q\rvert w_{0}v_{n-1}\ge w_{1}v_{n}>0$ for
  $n\ge1$. By Lemma~\ref{lem:positivity}, $A_{v}, A_{w}>0$, so
  $A_{v*w}=A_{v}A_{w}>0$ by Proposition~\ref{prop:prodbinet}.
  Now apply Lemma~\ref{lem:convergence}.
\end{proof}

\begin{remk}\label{remk:counterexamples}
  \begin{enumerate}
  \item Positive initial terms alone do not suffice. For
    $f(t)=t^{2}-3t+2$ (which satisfies $(C)$; $\theta_{1}=2$,
    $\theta_{2}=1$), $w_{0}=w_{1}=1$ gives $w_{n}=1$ for all $n$, and
    $\zeta_{w}(s)$ diverges for every $s$ (here $A=0$). The choice
    $w_{0}=2$, $w_{1}=1$ gives $w_{n}=3-2^{n}<0$ for $n\ge2$. For the unit
    case $f(t)=t^{2}-3t+1$, $w_{0}=3$, $w_{1}=1$ gives $w_{2}=0$.
  \item Proposition~\ref{prop:prodpositive} fails for $Q>0$. For
    $f(t)=t^{2}-3t+1$ ($\theta_{1}=\varphi^{2}$,
    $\varphi=(1+\sqrt5)/2$) and $v=w$ with $v_{0}=v_{1}=1$, $v$
    satisfies the hypotheses of Lemma~\ref{lem:positivity}, but
    $v*w$ has $(v*w)_{0}=-1$, $(v*w)_{1}=0$, so $\zeta_{v*w}(s)$ is not
    defined. Since $A_{v*w}=A_{v}A_{w}>0$, we still have $(v*w)_{n}>0$ for
    all sufficiently large $n$.
  \end{enumerate}
\end{remk}

\begin{assumption}\label{hyp:positivity}
  Throughout, the Lucas zeta function $\zeta_{u}(s)$ of a sequence
  $u\in\seqsp{f}{\realR}$ (in particular $u=v*w$ or $u=w^{*k}$) is
  considered only under the hypothesis
  \begin{equation}
    \label{eq:positivity}
    u_{n}>0\quad\text{for all $n\ge1$},\qquad A_{u}=u_{1}-u_{0}\theta_{2}>0,
  \end{equation}
  so that $\zeta_{u}(s)$ is defined and holomorphic on $\Re s>0$ by
  Lemma~\ref{lem:convergence}. If $Q<0$, then \eqref{eq:positivity} holds
  for $u=v*w$ whenever $v$ and $w$ satisfy the hypotheses of
  Proposition~\ref{prop:prodpositive}. 

  If $Q>0$, then
  \eqref{eq:positivity} for $v*w$ is an assumption; it does not follow
  from \eqref{eq:positivity} for $v$ and $w$
  (Remark~\ref{remk:counterexamples}~(2)). For the powers $u=w^{*k}$, see
  Proposition~\ref{prop:powerpositive}.
\end{assumption}

\begin{prop}\label{prop:powerpositive}
  Assume $(C)$, and let $w\in\seqsp{f}{\realR}$ with $A=A_{w}>0$. Put
  $B=B_{w}$ and $\lambda=B/A$.
  \begin{enumerate}
  \item If $\lvert\lambda\rvert\le1$, then $w^{*k}$ satisfies
    Hypothesis~\ref{hyp:positivity} for every $k\ge0$.
  \item If $\lvert\lambda\rvert>1$, then $w^{*k}$ does not satisfy
    Hypothesis~\ref{hyp:positivity} for any sufficiently large even $k$.
  \end{enumerate}
\end{prop}
\begin{proof}
  By Lemma~\ref{lem:kamano}, $d>0$ and $\theta_{1}>\lvert\theta_{2}\rvert>0$.
  By Proposition~\ref{prop:prodbinet}, for $k\ge0$ we have
  $A_{w^{*k}}=A^{k}>0$ and
  \[
    \sqrt d\,(w^{*k})_{n}=A^{k}\bigl(\theta_{1}^{n}-\lambda^{k}\theta_{2}^{n}\bigr).
  \]
  (1) For $n\ge1$, $\lvert\lambda^{k}\theta_{2}^{n}\rvert\le\lvert\theta_{2}\rvert^{n}<\theta_{1}^{n}$,
  so $(w^{*k})_{n}>0$.
  (2) Put $n=1$ if $Q>0$ and $n=2$ if $Q<0$. Since $\theta_{2}=Q/\theta_{1}$,
  we have $\theta_{2}^{n}>0$. For even $k$, $\lambda^{k}=\lvert\lambda\rvert^{k}$,
  so $(w^{*k})_{n}<0$ as soon as
  $\lvert\lambda\rvert^{k}>(\theta_{1}/\theta_{2})^{n}$.
\end{proof}

\subsection{Continuation, poles and residues}
\label{sec:lucas-zeta-functions-1}

Throughout this subsection we assume the condition $(C)$
(Definition~\ref{def:condc}), so that
$d>0$, $\theta_{1}>1$ and $\theta_{1}>\lvert\theta_{2}\rvert>0$
(Lemma~\ref{lem:kamano}). Put
\begin{equation}
  \label{eq:xm}
  r\coloneqq\frac{\lvert\theta_{2}\rvert}{\theta_{1}}\in(0,1),\qquad
  x_{m}(s)\coloneqq Q^{m}\theta_{1}^{-s-2m}\quad(m\in\intZ_{\ge0},\ s\in\cpxC),
\end{equation}
where $\theta_{1}^{-s}=\exp(-s\log\theta_{1})$. For $w\in\seqsp{f}{\realR}$
put $A_{w}=w_{1}-w_{0}\theta_{2}$, $B_{w}=w_{1}-w_{0}\theta_{1}$, and,
if $A_{w}\neq0$,
\begin{equation}
  \label{eq:lambda}
  \lambda_{w}\coloneqq\frac{B_{w}}{A_{w}}.
\end{equation}
For $\lambda\in\realR$ and an integer $n_{1}\ge1$ with
$\lvert\lambda\rvert r^{n_{1}}<1$, put
\begin{equation}
  \label{eq:Phi}
  \Phi_{n_{1}}(s;\lambda)\coloneqq\sum_{m=0}^{\infty}\binom{-s}{m}(-\lambda)^{m}
  \frac{x_{m}(s)^{n_{1}}}{1-x_{m}(s)}.
\end{equation}
Here $\binom{z}{m}=z(z-1)\cdots(z-m+1)/m!$. Note that
$\binom{-s}{m}(-1)^{m}=\binom{s+m-1}{m}$ and
$x_{m}(s)/(1-x_{m}(s))=Q^{m}/(\theta_{1}^{s+2m}-Q^{m})$. Finally let
\begin{equation}
  \label{eq:poleset}
  S\coloneqq\set{s_{m,j} | m\in\intZ_{\ge0},\ j\in\intZ},\qquad
  s_{m,j}\coloneqq-2m+\frac{m\log\lvert Q\rvert+(2j+l_{Q,m})\pi\sqrt{-1}}{\log\theta_{1}},
\end{equation}
with $l_{Q,m}=m$ if $Q<0$ and $l_{Q,m}=0$ if $Q>0$. Since
$Q^{m}=\lvert Q\rvert^{m}\exp(l_{Q,m}\pi\sqrt{-1})$, we have $x_{m}(s)=1$ if
and only if $(s+2m)\log\theta_{1}=m\log\lvert Q\rvert+(2j+l_{Q,m})\pi\sqrt{-1}$
for some $j\in\intZ$, that is, if and only if $s=s_{m,j}$.

\begin{lem}\label{lem:poleset}
  \begin{enumerate}
  \item $\Re s_{m,j}=-m\,\log(\theta_{1}/\lvert\theta_{2}\rvert)/\log\theta_{1}$.
    Hence $s_{m,j}=s_{m',j'}$ implies $(m,j)=(m',j')$, and $S$ is a
    discrete closed subset of $\cpxC$.
  \item $1-x_{m}(s)$ has a simple zero at each $s_{m,j}$, with derivative
    $\log\theta_{1}$ there.
  \end{enumerate}
\end{lem}
\begin{proof}
  (1) Since $\lvert Q\rvert=\theta_{1}\lvert\theta_{2}\rvert$,
  $\Re s_{m,j}=-2m+m(\log\theta_{1}+\log\lvert\theta_{2}\rvert)/\log\theta_{1}
  =-m\log(\theta_{1}/\lvert\theta_{2}\rvert)/\log\theta_{1}$, which is
  strictly decreasing in $m$ because $\theta_{1}>\lvert\theta_{2}\rvert$
  and $\theta_{1}>1$. For fixed $m$ the imaginary parts
  $(2j+l_{Q,m})\pi/\log\theta_{1}$ are distinct. A bounded region meets
  only finitely many $s_{m,j}$.
  
  (2) $\frac{\mathrm d}{\mathrm ds}(1-x_{m}(s))=x_{m}(s)\log\theta_{1}$, which equals
  $\log\theta_{1}\neq0$ where $x_{m}(s)=1$.
\end{proof}

\begin{thm}\label{thm:expforseq}
  Assume $(C)$, and let $w\in\seqsp{f}{\realR}$ satisfy
  Hypothesis~\ref{hyp:positivity}: $w_{n}>0$ for $n\ge1$ and
  $A=A_{w}>0$. Put $\lambda=\lambda_{w}$, and let $n_{1}\ge1$ be an integer
  with $\lvert\lambda\rvert r^{n_{1}}<1$.
  \begin{enumerate}
  \item For $\Re s>0$,
    \begin{equation}
      \label{eq:zetaw}
      \zeta_{w}(s)=\sum_{n=1}^{n_{1}-1}w_{n}^{-s}
      +\frac{d^{s/2}}{A^{s}}\,\Phi_{n_{1}}(s;\lambda),
    \end{equation}
    where $w_{n}^{-s}$, $d^{s/2}$, $A^{-s}$ are defined with real
    logarithms.
  \item The series $\Phi_{n_{1}}(s;\lambda)$ converges absolutely and
    uniformly on compact subsets of $\cpxC\setminus S$. Hence the
    right-hand side of \eqref{eq:zetaw} continues $\zeta_{w}(s)$ to a
    holomorphic function on $\cpxC\setminus S$. The continuation does not
    depend on $n_{1}$.
  \item The continuation is meromorphic on $\cpxC$. Every pole is simple
    and lies in $S$. At $s=s_{m,j}$,
    \[
      \mathop{\mathrm{Res}}_{s=s_{m,j}}\zeta_{w}(s)
      =\frac{d^{s_{m,j}/2}}{A^{s_{m,j}}}\binom{-s_{m,j}}{m}
      \frac{(-\lambda)^{m}}{\log\theta_{1}};
    \]
    thus $s_{m,j}$ is a pole if and only if
    $\binom{-s_{m,j}}{m}\lambda^{m}\neq0$.
  \end{enumerate}
  If $\lvert B\theta_{2}\rvert<A\theta_{1}$, one may take $n_{1}=1$, and
  \eqref{eq:zetaw} reads
  \[
    \zeta_{w}(s)=\frac{d^{s/2}}{A^{s}}\sum_{m=0}^{\infty}
    \binom{s+m-1}{m}\left(\frac{B}{A}\right)^{m}
    \frac{Q^{m}}{\theta_{1}^{s+2m}-Q^{m}}.
  \]
\end{thm}
\begin{proof}
  (1) Fix $s$ with $\tau=\Re s>0$, and put
  $\rho=\lvert\lambda\rvert r^{n_{1}}<1$. By Proposition~\ref{prop:binet},
  $w_{n}=A\theta_{1}^{n}(1-y_{n})/\sqrt{d}$ with
  $y_{n}=\lambda(\theta_{2}/\theta_{1})^{n}\in\realR$. For $n\ge n_{1}$
  we have $\lvert y_{n}\rvert\le\rho$, so $1-y_{n}>0$, and since
  $A, \theta_{1}, \sqrt d>0$,
  \[
    w_{n}^{-s}=\frac{d^{s/2}}{A^{s}}\,\theta_{1}^{-ns}(1-y_{n})^{-s}
    =\frac{d^{s/2}}{A^{s}}\,\theta_{1}^{-ns}
    \sum_{m=0}^{\infty}\binom{-s}{m}(-y_{n})^{m},
  \]
  by the binomial series, valid for real $\lvert y_{n}\rvert<1$ with the real
  logarithm of $1-y_{n}$. Since
  $\bigl\lvert\binom{-s}{m}\bigr\rvert
  =\lvert s(s+1)\cdots(s+m-1)\rvert/m!
  \le\binom{\lvert s\rvert+m-1}{m}$, we get
  \begin{align*}
    \sum_{n\ge n_{1}}\sum_{m\ge0}
    \Bigl\lvert\theta_{1}^{-ns}\binom{-s}{m}(-y_{n})^{m}\Bigr\rvert
    &\le\sum_{n\ge n_{1}}\theta_{1}^{-n\tau}
      \sum_{m\ge0}\binom{\lvert s\rvert+m-1}{m}\rho^{m}\\
    &=\frac{(1-\rho)^{-\lvert s\rvert}\,\theta_{1}^{-n_{1}\tau}}{1-\theta_{1}^{-\tau}}<\infty.
  \end{align*}
  So the double series converges absolutely, and we may sum over $n$
  first. As $\theta_{2}=Q/\theta_{1}$,
  $\theta_{1}^{-s}(\theta_{2}/\theta_{1})^{m}=x_{m}(s)$ and
  $\lvert x_{m}(s)\rvert=\theta_{1}^{-\tau}r^{m}<1$, so
  $\sum_{n\ge n_{1}}\theta_{1}^{-ns}(-y_{n})^{m}
  =(-\lambda)^{m}x_{m}(s)^{n_{1}}/(1-x_{m}(s))$. Adding the terms with
  $n<n_{1}$ gives \eqref{eq:zetaw}.

  (2) Let $\Omega\subset\cpxC\setminus S$ be compact, and put
  $\kappa=\max_{\Omega}\lvert s\rvert$ and $\tau_{0}=\min_{\Omega}\Re s$. On $\Omega$,
  $\lvert x_{m}(s)\rvert\le\theta_{1}^{-\tau_{0}}r^{m}$, so there is $m_{0}$ with
  $\lvert x_{m}(s)\rvert\le1/2$ on $\Omega$ for $m\ge m_{0}$. For $m\ge m_{0}$ and $s\in\Omega$
  the $m$-th term of \eqref{eq:Phi} is bounded by
  $2\theta_{1}^{-n_{1}\tau_{0}}\binom{\kappa+m-1}{m}\rho^{m}$, and
  $\sum_{m}\binom{\kappa+m-1}{m}\rho^{m}=(1-\rho)^{-\kappa}<\infty$. The
  finitely many terms with $m<m_{0}$ are holomorphic on a neighborhood of
  $\Omega$, because $1-x_{m}(s)$ vanishes only on $S$. By the Weierstrass
  theorem, $\Phi_{n_{1}}(s;\lambda)$ is holomorphic on $\cpxC\setminus S$.
  The other two parts of the right-hand side of \eqref{eq:zetaw} are entire,
  and the right-hand side agrees with $\zeta_{w}(s)$ on $\Re s>0$ by (1). By
  the identity theorem, the continuation is unique and therefore
  independent of $n_{1}$.

  (3) Fix $s^{*}=s_{m,j}\in S$ and a closed disc $U$ around $s^{*}$
  with $U\cap S=\set{s^{*}}$ (Lemma~\ref{lem:poleset}~(1)). For
  $m'\neq m$, $1-x_{m'}$ has no zero on $U$. The estimate in (2)
  applies on $U$ to $\sum_{m'\neq m}$, so this sum is holomorphic on
  $U$. The $m$-th term is $\psi(s)/(1-x_{m}(s))$ with
  $\psi(s)=\binom{-s}{m}(-\lambda)^{m}x_{m}(s)^{n_{1}}$ entire, and
  $x_{m}(s^{*})=1$. By Lemma~\ref{lem:poleset}~(2) its residue is
  $\psi(s^{*})/\log\theta_{1}=\binom{-s^{*}}{m}(-\lambda)^{m}/\log\theta_{1}$.
  Multiplying by the entire, nowhere vanishing factor $d^{s/2}A^{-s}$ gives
  the residue. So the continuation is meromorphic at every point of $S$.
  The last statement is (1) with $n_{1}=1$, since
  $\lvert\lambda\rvert r<1$ means $\lvert B\theta_{2}\rvert<A\theta_{1}$.
\end{proof}

\begin{remk}
  For $w=\calF$, we have $A=B=1$ and $\lambda=1$, and
  Theorem~\ref{thm:expforseq} with $n_{1}=1$ is Kamano's
  Theorem~4~\cite{MR3064567}. The expansion by the binomial series is
  the standard method for these series; see
  also~\cite{MR4594520,AKLW2024}. In general one cannot take $n_{1}=1$. For example, let $f(t)=t^{2}-3t+1$ and $(w_{0},w_{1})=(5,2)$;
  this $w$ satisfies Hypothesis~\ref{hyp:positivity} by
  Lemma~\ref{lem:positivity}. With $\varphi=(1+\sqrt5)/2$ we have
  $\theta_{1}=\varphi^{2}$, $\theta_{2}=\varphi^{-2}$, $r=\varphi^{-4}$,
  $A=(5\sqrt5-11)/2=\varphi^{-5}$ and $B=-(5\sqrt5+11)/2=-\varphi^{5}$.
  Hence $\lvert\lambda\rvert r^{n}=\varphi^{10-4n}$, and $n_{1}=3$ is the
  least admissible value.
\end{remk}

\begin{prop}\label{prop:lambda}
  Assume $d>0$. The map $\realR[t]/(f(t))\to\realR\times\realR$,
  $w_{1}-w_{0}t\mapsto(A_{w},B_{w})$, is an isomorphism of $\realR$-algebras.
  Consequently:
  \begin{enumerate}
  \item $A_{v*w}=A_{v}A_{w}$ and $B_{v*w}=B_{v}B_{w}$ for all
    $v, w\in\seqsp{f}{\realR}$.
  \item $\bvec{w_{1}}{w_{0}}\mapsto\lambda_{w}=B_{w}/A_{w}$~\eqref{eq:lambda} is a well defined
    group isomorphism $G_{\realR}(f)\to\realR^{\times}$.
  \end{enumerate}
  Part~(2) is the case $R=\realR$ of the isomorphism of Aoki and
  Kida~\cite{MR3902117} for local rings $R$ with $d\in R^{\times}$.
\end{prop}
\begin{proof}
  The map is evaluation at the two distinct real roots $t=\theta_{2}$ and
  $t=\theta_{1}$ of $f$. It is an isomorphism by the Chinese remainder
  theorem, since $(t-\theta_{1})$ and $(t-\theta_{2})$ are coprime. Part (1)
  is Proposition~\ref{prop:prodbinet}. For (2): the units
  of $\realR\times\realR$ are the pairs $(A,B)$ with $AB\neq0$. Scaling $w$
  by $c\in\realR^{\times}$ scales $(A,B)$ by $c$, so it does not change
  $B/A$. The induced map $(\realR^{\times})^{2}/\realR^{\times}\to\realR^{\times}$,
  $(A,B)\mapsto B/A$, is an isomorphism, and it is a homomorphism by (1).
\end{proof}

\begin{cor}\label{cor:zetaprod}
  Assume $(C)$. Let $v, w\in\seqsp{f}{\realR}$, and assume that $v$, $w$
  and $u=v*w$ satisfy Hypothesis~\ref{hyp:positivity}. (If $Q<0$,
  $v_{0},w_{0}\ge0$ and $v_{1},w_{1}>0$, then the hypothesis on $u$ holds by
  Proposition~\ref{prop:prodpositive}.)
  Let $n_{1}\ge1$ satisfy
  $\lvert\lambda_{v}\lambda_{w}\rvert r^{n_{1}}<1$, with $\lambda$ as
  in~\eqref{eq:lambda} and $r$ as in~\eqref{eq:xm}. Then, with
  $\Phi_{n_{1}}$ as in~\eqref{eq:Phi}, as meromorphic
  functions on $\cpxC$,
  \begin{equation}
    \label{eq:zetavw}
    \zeta_{v*w}(s)=\sum_{n=1}^{n_{1}-1}u_{n}^{-s}
    +\frac{d^{s/2}}{(A_{v}A_{w})^{s}}\,\Phi_{n_{1}}(s;\lambda_{v}\lambda_{w}).
  \end{equation}
  The poles of $\zeta_{v*w}$ lie in the same set $S$~\eqref{eq:poleset} as those of $\zeta_{v}$
  and $\zeta_{w}$, and
  \[
    \mathop{\mathrm{Res}}_{s=s_{m,j}}\zeta_{v*w}(s)
    =\frac{d^{s_{m,j}/2}}{(A_{v}A_{w})^{s_{m,j}}}\binom{-s_{m,j}}{m}
    \frac{(-\lambda_{v}\lambda_{w})^{m}}{\log\theta_{1}}.
  \]
\end{cor}
\begin{proof}
  Apply Theorem~\ref{thm:expforseq} to $u$, and use
  $A_{u}=A_{v}A_{w}$ and $\lambda_{u}=\lambda_{v}\lambda_{w}$
  (Proposition~\ref{prop:lambda}).
\end{proof}

\begin{remk}\label{remk:zetaprodstructure}
  Theorem~\ref{thm:expforseq} and Corollary~\ref{cor:zetaprod} show that
  $\zeta_{w}$ depends on the class of $w$ only through $A_{w}$ and
  $\lambda_{w}$.
  \begin{itemize}
  \item Up to the elementary factor $d^{s/2}A^{-s}$ and finitely many
    initial terms, all three are values of the same function
    $\Phi_{n_{1}}(s;\lambda)$. It is evaluated at $\lambda_{v}$,
    $\lambda_{w}$ and $\lambda_{v}\lambda_{w}$, where
    $\lambda\colon\GR{f}\to\realR^{\times}$ is the isomorphism of
    Proposition~\ref{prop:lambda}.
  \item As a power series in $\lambda$,
    $\Phi_{n_{1}}(s;\lambda)=\sum_{m}a_{m}(s)\lambda^{m}$, convergent for
    $\lvert\lambda\rvert<r^{-n_{1}}$, where
    \[
      a_{m}(s)=\binom{-s}{m}(-1)^{m}\frac{x_{m}(s)^{n_{1}}}{1-x_{m}(s)}.
    \] Passing from $v$ to $v*w$
    replaces $\lambda_{v}^{m}$ by $(\lambda_{v}\lambda_{w})^{m}$ in the
    $m$-th term, and $A_{v}^{-s}$ by $(A_{v}A_{w})^{-s}$.
  \item In particular, if $n_{1}=1$ is admissible for $w$, the normalized
    function $\widehat\zeta_{w}(s)\coloneqq(A_{w}/\sqrt{d})^{s}\zeta_{w}(s)$
    equals $\Phi_{1}(s;\lambda_{w})$. Hence, if $n_{1}=1$ is admissible
    for $v*w$, then $\widehat\zeta_{v*w}(s)=\Phi_{1}(s;\lambda_{v}\lambda_{w})$.
  \end{itemize}
\end{remk}

\begin{cor}\label{cor:zetapower}
  Assume $(C)$, let $k\in\intZ$, and let $w\in\seqsp{f}{\realR}$, with
  $A_{w}B_{w}\neq0$ if $k<0$, be such that $w^{*k}$ satisfies
  Hypothesis~\ref{hyp:positivity}. Let $n_{1}\ge1$
  satisfy $\lvert\lambda_{w}\rvert^{k}r^{n_{1}}<1$ (notation
  of~\eqref{eq:xm}, \eqref{eq:lambda} and~\eqref{eq:Phi}). Then
  \begin{equation}
    \label{eq:zetawkth}
    \zeta_{w^{*k}}(s)=\sum_{n=1}^{n_{1}-1}(w^{*k})_{n}^{-s}
    +\frac{d^{s/2}}{(A_{w}^{k})^{s}}\,\Phi_{n_{1}}(s;\lambda_{w}^{k}).
  \end{equation}
\end{cor}
\begin{proof}
  By Proposition~\ref{prop:lambda}, $A_{w^{*k}}=A_{w}^{k}$ and
  $\lambda_{w^{*k}}=\lambda_{w}^{k}$. Apply Theorem~\ref{thm:expforseq}.
\end{proof}

\begin{ex}
  For $k=2$ and $n_{1}=1$ (admissible if $\lambda_{w}^{2}r<1$),
  \begin{equation}
    \label{eq:zetaww}
    \zeta_{w^{*2}}(s)
    =\frac{d^{s/2}}{A^{2s}}\sum_{m=0}^{\infty}
    \binom{-s}{m}(-1)^{m}\left(\frac{B}{A}\right)^{2m}
    \frac{Q^{m}}{\theta_{1}^{s+2m}-Q^{m}}.
  \end{equation}
  Take the companion Lucas sequence $\calL$, with
  $\pvec{\calL_{1}}{\calL_{0}}=\pvec{P}{2}$. Then
  $A=P-2\theta_{2}=\sqrt{d}$ and $B=P-2\theta_{1}=-\sqrt{d}$, so
  $\lambda_{\calL}=-1$ and $\lambda_{\calL}^{2}r=r<1$. We also have
  $\calL^{*2}=d\calF$, because $A_{\calL}^{2}=B_{\calL}^{2}=d$ and
  $\calF$ has $A=B=1$. By \eqref{eq:zetaww},
  \[
    \zeta_{\calL^{*2}}(s)
    =d^{-s/2}\sum_{m=0}^{\infty}\binom{-s}{m}(-1)^{m}
    \frac{Q^{m}}{\theta_{1}^{s+2m}-Q^{m}}=d^{-s}\zeta_{\calF}(s),
  \]
  as Remark~\ref{remk:scaling} predicts.
\end{ex}

\begin{remk}
  The proof of Theorem~\ref{thm:expforseq} parallels a computation of
  Kaneko, Kurokawa and Wakayama~\cite{MR2069738}: they expand
  $(1-q^{n})^{-s}$ in the same way for a $q$-analog $\zeta_{q}(s)$ of the
  Riemann zeta function $\zeta(s)$ ($0<q<1$), and prove that
  $\zeta_{q}(s)\to\zeta(s)$ as $q$ tends to $1$ from below, for
  $s\neq1$~\cite[Theorem~2]{MR2069738}.
\end{remk}

\subsection{Partial Lucas zeta functions and Lucas \texorpdfstring{$L$}{L}-functions}
\label{sec:lucas-l-functions}

We keep the notation and the standing assumption $(C)$ of
\S\ref{sec:lucas-zeta-functions-1}: in particular $r$ and $x_{m}(s)$ are
defined in~\eqref{eq:xm}, $\lambda_{w}$ in~\eqref{eq:lambda},
$\Phi_{n_{1}}$ in~\eqref{eq:Phi}, and $l_{Q,m}$ and the pole set $S$
in~\eqref{eq:poleset}. Throughout, $N\ge1$ is an integer
and $\xi_{N}=\exp(2\pi\sqrt{-1}/N)$. This subsection extends
Kamano~\cite[\S4]{MR3064567}, where $w=\calF$, to general $w$.

\begin{defn}\label{def:partial}
  For $w\in\seqsp{f}{\realR}$ satisfying Hypothesis~\ref{hyp:positivity}
  and $a\in\intZ$, the \emph{partial Lucas zeta function} is
  \[
    \zeta_{w}(s,a,N)\coloneqq\sum_{\substack{n\ge1\\ n\equiv a\pmod{N}}}w_{n}^{-s}
    \qquad(\Re s>0).
  \]
  It is a subseries of $\zeta_{w}(s)$, so by Lemma~\ref{lem:convergence}
  it converges absolutely and locally uniformly on $\Re s>0$.
\end{defn}

Put
\[
  s^{(N)}_{m,j}\coloneqq-2m+\frac{m\log\lvert Q\rvert+(2j/N+l_{Q,m})\pi\sqrt{-1}}{\log\theta_{1}}
  \quad(m\in\intZ_{\ge0},\ j\in\intZ),
\]
and let $S_{N}$ be the set of these points; thus $S_{1}=S$.

\begin{lem}\label{lem:polesetN}
  \begin{enumerate}
  \item $x_{m}(s)^{N}=1$ if and only if $s=s^{(N)}_{m,j}$ for some
    $j\in\intZ$, and then $x_{m}(s)=\xi_{N}^{-j}$.
  \item $s^{(N)}_{m,j}=s^{(N)}_{m',j'}$ implies $(m,j)=(m',j')$, and
    $S_{N}$ is discrete and closed in $\cpxC$.
  \item $1-x_{m}(s)^{N}$ has a simple zero at each $s^{(N)}_{m,j}$, with
    derivative $N\log\theta_{1}$ there.
  \end{enumerate}
\end{lem}
\begin{proof}
  (1) Since $Q^{m}=\lvert Q\rvert^{m}\exp(l_{Q,m}\pi\sqrt{-1})$, we have
  $x_{m}(s)=\exp(\Lambda_{m}(s))$ with
  \[
    \Lambda_{m}(s)=m\log\lvert Q\rvert+l_{Q,m}\pi\sqrt{-1}-(s+2m)\log\theta_{1}.
  \]
  Hence $x_{m}(s)^{N}=1$ if and only if $N \Lambda_{m}(s)=2\pi\sqrt{-1}\,j$ for
  some $j\in\intZ$. Solving for $s$ gives $s=s^{(N)}_{m,-j}$, and then
  $x_{m}(s)=\exp(2\pi\sqrt{-1}j/N)$.

  (2) $\Re s^{(N)}_{m,j}=\Re s_{m,0}$ is strictly decreasing in $m$
  (Lemma~\ref{lem:poleset}~(1)), and for fixed $m$ the imaginary parts
  are distinct. A bounded set meets finitely many points.

  (3) $\frac{\mathrm d}{\mathrm ds}(1-x_{m}(s)^{N})=N x_{m}(s)^{N}\log\theta_{1}$, which equals
  $N\log\theta_{1}$ where $x_{m}(s)^{N}=1$.
\end{proof}

\begin{thm}\label{thm:partial}
  Assume $(C)$, and let $w$ satisfy Hypothesis~\ref{hyp:positivity}. Put
  $A=A_{w}$ and $\lambda=\lambda_{w}$. Let $n_{1}\ge1$ satisfy
  $\lvert\lambda\rvert r^{n_{1}}<1$, let $a\in\intZ$, and let $a'$ be the
  least integer $\ge n_{1}$ with $a'\equiv a\pmod N$.
  \begin{enumerate}
  \item For $\Re s>0$,
    \[
      \zeta_{w}(s,a,N)=\sum_{\substack{1\le n<n_{1}\\ n\equiv a\ (\mathrm{mod}\ N)}}w_{n}^{-s}
      +\frac{d^{s/2}}{A^{s}}\sum_{m=0}^{\infty}\binom{-s}{m}(-\lambda)^{m}
      \frac{x_{m}(s)^{a'}}{1-x_{m}(s)^{N}}.
    \]
  \item The right-hand side converges absolutely and uniformly on compact
    subsets of $\cpxC\setminus S_{N}$, and continues $\zeta_{w}(s,a,N)$ to a
    meromorphic function on $\cpxC$ that is holomorphic outside $S_{N}$.
  \item Every pole is simple and lies in $S_{N}$. At $s^{*}=s^{(N)}_{m,j}$,
    \begin{equation}
      \label{eq:respartial}
      \mathop{\mathrm{Res}}_{s=s^{*}}\zeta_{w}(s,a,N)
      =\frac{d^{s^{*}/2}}{A^{s^{*}}}\binom{-s^{*}}{m}(-\lambda)^{m}
      \frac{\xi_{N}^{-ja}}{N\log\theta_{1}}.
    \end{equation}
  \end{enumerate}
\end{thm}
\begin{proof}
  The proof of Theorem~\ref{thm:expforseq} applies with the following
  changes.

  In (1), the sum over $n\ge n_{1}$ is restricted to
  $n=a'+N\nu$, $\nu\ge0$. The absolute bound used there dominates this
  subseries. The geometric sum becomes
  $\sum_{\nu\ge0}x_{m}(s)^{a'+N\nu}=x_{m}(s)^{a'}/(1-x_{m}(s)^{N})$, valid since
  $\lvert x_{m}(s)\rvert<1$ for $\Re s>0$.

  In (2), on a compact set $\Omega\subset\cpxC\setminus S_{N}$ with
  $\tau_{0}=\min_{\Omega}\Re s$ and $\kappa=\max_{\Omega}\lvert s\rvert$, we have
  $\lvert x_{m}(s)^{a'}\rvert\le\theta_{1}^{a'\lvert\tau_{0}\rvert}r^{ma'}\le\theta_{1}^{a'\lvert\tau_{0}\rvert}r^{mn_{1}}$.
  Also $\lvert x_{m}(s)^{N}\rvert\le1/2$ on $\Omega$ for $m\ge m_{0}$. So the
  $m$-th term is bounded by
  $2\theta_{1}^{a'\lvert\tau_{0}\rvert}\binom{\kappa+m-1}{m}\rho^{m}$ with
  $\rho=\lvert\lambda\rvert r^{n_{1}}<1$. The terms with $m<m_{0}$ are
  holomorphic on a neighborhood of $\Omega$ by Lemma~\ref{lem:polesetN}~(1).

  In (3), use Lemma~\ref{lem:polesetN}~(2), (3) in place of
  Lemma~\ref{lem:poleset}. The numerator at $s^{*}$ is
  $\binom{-s^{*}}{m}(-\lambda)^{m}x_{m}(s^{*})^{a'}$, and
  $x_{m}(s^{*})^{a'}=\xi_{N}^{-ja'}=\xi_{N}^{-ja}$ by
  Lemma~\ref{lem:polesetN}~(1), because $a'\equiv a\pmod N$.
\end{proof}

\begin{defn}\label{def:Lfunction}
  Let $\chi$ be a Dirichlet character modulo $N$, and let $w$ satisfy
  Hypothesis~\ref{hyp:positivity}. The \emph{Lucas $L$-function} is
  \[
    L_{w}(s,\chi)\coloneqq\sum_{n=1}^{\infty}\frac{\chi(n)}{w_{n}^{s}}
    =\sum_{a=1}^{N}\chi(a)\,\zeta_{w}(s,a,N)\qquad(\Re s>0).
  \]
  For $j\in\intZ$ put $\mathcal G_{j}(\chi)\coloneqq\sum_{a=1}^{N}\chi(a)\xi_{N}^{-ja}$.
\end{defn}

\begin{thm}\label{thm:resL}
  Under the assumptions of Theorem~\ref{thm:partial}, $L_{w}(s,\chi)$
  continues to a meromorphic function on $\cpxC$. Its poles are simple and
  lie in $S_{N}$, and
  \[
    \mathop{\mathrm{Res}}_{s=s^{(N)}_{m,j}}L_{w}(s,\chi)
    =\frac{d^{s^{*}/2}}{A^{s^{*}}}\binom{-s^{*}}{m}(-\lambda)^{m}
    \frac{\mathcal G_{j}(\chi)}{N\log\theta_{1}}\qquad(s^{*}=s^{(N)}_{m,j}).
  \]
\end{thm}
\begin{proof}
  $L_{w}$ is a finite linear combination of the $\zeta_{w}(s,a,N)$. Apply
  Theorem~\ref{thm:partial} and sum \eqref{eq:respartial} with the
  weights $\chi(a)$.
\end{proof}

\begin{cor}\label{cor:Lreal}
  If $\chi$ is not the principal character (the character with
  $\chi(a)=1$ for all $a$ prime to $N$), then $L_{w}(s,\chi)$ is
  holomorphic at every real $s$.
\end{cor}
\begin{proof}
  By Theorem~\ref{thm:resL}, it suffices to show $\mathcal G_{j}(\chi)=0$ whenever
  $s^{(N)}_{m,j}$ is real. The point is real exactly when
  $2j/N+l_{Q,m}=0$, i.e.\ $2j=-Nl_{Q,m}$.
  \begin{itemize}
  \item If $l_{Q,m}$ is even, then $N\mid j$, so $\xi_{N}^{-ja}=1$ and
    $\mathcal G_{j}(\chi)=\sum_{a=1}^{N}\chi(a)=0$ by the orthogonality of
    characters.
  \item If $l_{Q,m}$ is odd, then $N$ is even and $j/N=-l_{Q,m}/2$. Hence
    $\xi_{N}^{-ja}=\exp(\pi\sqrt{-1}\,l_{Q,m}a)$, which is $(-1)^{a}$. Since $N$ is
    even, $\chi(a)=0$ for even $a$. Hence
    $\mathcal G_{j}(\chi)=-\sum_{a=1}^{N}\chi(a)=0$.
  \end{itemize}
  For $w=\calF$ this is Kamano~\cite[Corollary~11]{MR3064567}.
\end{proof}

\begin{cor}\label{cor:Lprod}
  Let $v, w$, and $u=v*w$ satisfy Hypothesis~\ref{hyp:positivity}.
  Then Theorems~\ref{thm:partial} and~\ref{thm:resL} hold for $u$ with
  $A$ replaced by $A_{v}A_{w}$ and $\lambda$ by $\lambda_{v}\lambda_{w}$,
  where $n_{1}$ satisfies $\lvert\lambda_{v}\lambda_{w}\rvert r^{n_{1}}<1$.
  In particular,
  \[
    \mathop{\mathrm{Res}}_{s=s^{(N)}_{m,j}}L_{v*w}(s,\chi)
    =\frac{d^{s^{*}/2}}{(A_{v}A_{w})^{s^{*}}}\binom{-s^{*}}{m}
    (-\lambda_{v}\lambda_{w})^{m}\frac{\mathcal G_{j}(\chi)}{N\log\theta_{1}}.
  \]
  The same holds for $w^{*k}$ ($k\in\intZ$, with $A_{w}B_{w}\neq0$ if
  $k<0$) with $A_{w}^{k}$ and $\lambda_{w}^{k}$, when $w^{*k}$ satisfies
  Hypothesis~\ref{hyp:positivity}.
\end{cor}
\begin{proof}
  By Proposition~\ref{prop:lambda}, $A_{v*w}=A_{v}A_{w}$,
  $\lambda_{v*w}=\lambda_{v}\lambda_{w}$, $A_{w^{*k}}=A_{w}^{k}$ and
  $\lambda_{w^{*k}}=\lambda_{w}^{k}$.
\end{proof}

\section{Special values of Lucas zeta functions and Laxton groups}
\label{sec:special-values-lucas}

\subsection{Values at negative integers}
\label{sec:special-values-lucas-1}

We keep the assumption $(C)$ (Definition~\ref{def:condc}) and the
notation of \S\ref{sec:lucas-zeta-functions-1}: $\theta_{1}$,
$\theta_{2}$, $A_{w}$, $B_{w}$ (Proposition~\ref{prop:binet}), $x_{m}(s)$
and $r$~\eqref{eq:xm}, $\lambda_{w}$~\eqref{eq:lambda},
$\Phi_{n_{1}}$~\eqref{eq:Phi}, and $S$, $s_{m,j}$~\eqref{eq:poleset}. Let
$m\ge1$ be an integer, and put
\begin{equation}
  \label{eq:beta}
  \beta_{i}\coloneqq\theta_{1}^{m-i}\theta_{2}^{i}=x_{i}(-m)\qquad(i\in\intZ_{\ge0}),
\end{equation}
where the second equality holds because $Q^{i}=\theta_{1}^{i}\theta_{2}^{i}$. For
$w\in\seqsp{f}{\realR}$ and $n_{1}\ge1$ put
\begin{equation}
  \label{eq:Zmn1}
  Z_{w}(m;n_{1})\coloneqq\sum_{n=1}^{n_{1}-1}w_{n}^{m}
  +d^{-m/2}\sum_{i=0}^{m}\binom{m}{i}(-1)^{i}A_{w}^{m-i}B_{w}^{i}
  \frac{\beta_{i}^{n_{1}}}{1-\beta_{i}},
\end{equation}
whenever $\beta_{i}\neq1$ for $0\le i\le m$.

\begin{lem}\label{lem:alpha}
  \begin{enumerate}
  \item $-m\in S$ if and only if $\beta_{i}=1$ for some $i\ge0$; there is
    at most one such $i$.
  \item If $\beta_{i}=1$ with $i>m$, then $\lvert\theta_{2}\rvert>1$.
    Hence, if $\lvert\theta_{2}\rvert\le1$, then $\beta_{i}\neq1$ for all
    $i>m$.
  \item If $\beta_{i}\neq1$ for $0\le i\le m$, then $Z_{w}(m;n_{1})$ of~\eqref{eq:Zmn1} does
    not depend on $n_{1}\ge1$. We write $Z_{w}(m)$ for it.
  \item Let $P, Q\in\ratQ$ with $d$ not a square in $\ratQ$. If
    $\beta_{i}=1$, then $Q=\pm1$ and $i=m/2$. In particular,
    $\beta_{i}\neq1$ for all $i>m$; and $\beta_{i}\neq1$ for all
    $0\le i\le m$ unless $Q=1$ and $m$ is even, or $Q=-1$ and $4\mid m$.
  \end{enumerate}
\end{lem}
\begin{proof}
  (1) By the definition of $S$, $x_{i}(s)=1$ if and only if $s=s_{i,j}$
  for some $j$. Hence $-m\in S$ if and only if
  $\beta_{i}=x_{i}(-m)=1$ for some $i$. By Lemma~\ref{lem:poleset}~(1)
  the pair $(i,j)$ with $-m=s_{i,j}$ is unique.
  (2) $\lvert\beta_{i}\rvert=1$ gives
  $\lvert\theta_{2}\rvert^{i}=\theta_{1}^{i-m}>1$ for $i>m$, since
  $\theta_{1}>1$.
  (3) Write $D=Z_{w}(m;n_{1}+1)-Z_{w}(m;n_{1})$. Since
  $\beta_{i}^{n_{1}+1}-\beta_{i}^{n_{1}}=-(1-\beta_{i})\beta_{i}^{n_{1}}$,
  \[
    D=w_{n_{1}}^{m}-d^{-m/2}\sum_{i=0}^{m}\binom{m}{i}(-1)^{i}
    (A\theta_{1}^{n_{1}})^{m-i}(B\theta_{2}^{n_{1}})^{i}
    =w_{n_{1}}^{m}-\Bigl(\frac{A\theta_{1}^{n_{1}}-B\theta_{2}^{n_{1}}}{\sqrt{d}}\Bigr)^{m}=0
  \]
  by the binomial theorem and Proposition~\ref{prop:binet}.
  (4) Let $\sigma$ be the nontrivial automorphism of $\ratQ(\sqrt d)$; it
  exchanges $\theta_{1}$ and $\theta_{2}$. If $\beta_{i}=1$, then also
  $\sigma(\beta_{i})=\theta_{2}^{m-i}\theta_{1}^{i}=1$, and the product
  gives $(\theta_{1}\theta_{2})^{m}=Q^{m}=1$. Hence $Q=\pm1$ and
  $\lvert\theta_{2}\rvert=\theta_{1}^{-1}$, so
  $1=\lvert\beta_{i}\rvert=\theta_{1}^{m-2i}$ and $i=m/2$. Then
  $\beta_{m/2}=Q^{m/2}$, which equals $1$ exactly in the two cases stated.
\end{proof}

\begin{thm}\label{thm:negval}
  Assume $(C)$, let $w$ satisfy Hypothesis~\ref{hyp:positivity}, and let
  $m\ge1$. Put $A=A_{w}$, $B=B_{w}$, $\lambda=B/A$, and let $\beta_{i}$ be as
  in~\eqref{eq:beta}.
  \begin{enumerate}
  \item If $\beta_{i}=1$ for some $0\le i\le m$, then $-m=s_{i,j}$~\eqref{eq:poleset} for some
    $j$. The function $\zeta_{w}$ has at most a simple pole there, with the
    residue of Theorem~\ref{thm:expforseq}~(3).
  \item Otherwise $\zeta_{w}$ is holomorphic at $s=-m$, and
    \[
      \zeta_{w}(-m)=Z_{w}(m)+E_{w}(m),
    \]
    where $Z_{w}(m)$ is defined in Lemma~\ref{lem:alpha}~(3).
    Here $E_{w}(m)=0$ if $\beta_{i}\neq1$ for all $i>m$. If $\beta_{i}=1$
    for the (unique) $i>m$, then
    \[
      E_{w}(m)=\Bigl(\frac{A}{\sqrt{d}}\Bigr)^{m}(-1)^{i-m}
      \frac{m!\,(i-m-1)!}{i!}\,\frac{(-\lambda)^{i}}{\log\theta_{1}}.
    \]
  \end{enumerate}
  In particular, $\zeta_{w}(-m)=Z_{w}(m)$ if $\beta_{i}\neq1$ for
  $0\le i\le m$ and either $\lvert\theta_{2}\rvert\le1$, or
  $P, Q\in\ratQ$ and $d$ is not a square.
\end{thm}
\begin{proof}
  (1) is Lemma~\ref{lem:alpha}~(1) and Theorem~\ref{thm:expforseq}~(3).

  (2) Choose $n_{1}$ as in Theorem~\ref{thm:expforseq}. We use the
  expression \eqref{eq:zetaw}, i.e.\ the finite part plus
  $d^{s/2}A^{-s}\Phi_{n_{1}}(s;\lambda)$.

  First suppose $-m\notin S$. The series $\Phi_{n_{1}}$ converges
  uniformly on a closed disc around $-m$ disjoint from $S$
  (Theorem~\ref{thm:expforseq}~(2)), so its value at $-m$ is the sum of
  the values of its terms. At $s=-m$ we have $\binom{m}{i'}=0$ for $i'>m$
  and $x_{i'}(-m)=\beta_{i'}$. Also
  \[
    d^{-m/2}A^{m}\binom{m}{i'}(-\lambda)^{i'}
    =d^{-m/2}\binom{m}{i'}(-1)^{i'}A^{m-i'}B^{i'}.
  \]
  This gives $\zeta_{w}(-m)=Z_{w}(m;n_{1})=Z_{w}(m)$.

  Next suppose $-m\in S$. By (1) and Lemma~\ref{lem:alpha}~(1),
  $-m=s_{i,j}$ for a unique pair with $i>m$. As in the proof of
  Theorem~\ref{thm:expforseq}~(3), the series without its $i$-th term
  converges uniformly on a closed disc $U$ around $-m$ with
  $U\cap S=\set{-m}$. So that part contributes $Z_{w}(m;n_{1})$ as
  before. Its terms $i'>m$, $i'\neq i$, vanish at $-m$. The
  $i$-th term is $\psi(s)/(1-x_{i}(s))$ with
  $\psi(s)=\binom{-s}{i}(-\lambda)^{i}x_{i}(s)^{n_{1}}$.

  Write $\binom{-s}{i}=\frac{1}{i!}\prod_{h=0}^{i-1}(-s-h)$. Exactly one
  factor, $h=m$, vanishes at $s=-m$; it equals $-(s+m)$. The other factors
  give
  \[
    \prod_{h\ne m}(m-h)=m!\,(-1)^{i-1-m}(i-m-1)!.
  \]
  Hence
  \[
    \binom{-s}{i}=(s+m)\,(-1)^{i-m}\frac{m!\,(i-m-1)!}{i!}+O((s+m)^{2}).
  \]
  By Lemma~\ref{lem:poleset}~(2),
  $1-x_{i}(s)=(s+m)\log\theta_{1}+O((s+m)^{2})$, and $x_{i}(-m)=1$.
  Therefore $\psi(s)/(1-x_{i}(s))$ extends holomorphically to $-m$, with value
  \[
    (-1)^{i-m}\frac{m!\,(i-m-1)!}{i!}\,\frac{(-\lambda)^{i}}{\log\theta_{1}}.
  \]
  Multiplying by $d^{-m/2}A^{m}=(A/\sqrt{d})^{m}$ gives $E_{w}(m)$.

  The last assertion follows from Lemma~\ref{lem:alpha}~(2) and~(4).
\end{proof}

\begin{cor}\label{cor:negform}
  For $m\ge1$ with $\beta_{i}\neq1$ for $0\le i\le m$, put
  \[
    c_{i}(m)\coloneqq d^{-m/2}\binom{m}{i}(-1)^{i}\frac{\beta_{i}}{1-\beta_{i}}
    \quad(0\le i\le m),
  \]
  \[
    H_{m}(X,Y)\coloneqq\sum_{i=0}^{m}c_{i}(m)X^{m-i}Y^{i}.
  \]
  The binary form $H_{m}$ depends only on $f$ and $m$. Then
  $Z_{w}(m)=H_{m}(A_{w},B_{w})$. Moreover, for
  $v, w\in\seqsp{f}{\realR}$ with $A_{v}B_{v}A_{w}B_{w}\neq0$ and $\ell\in\intZ$,
  \[
    Z_{v*w}(m)=H_{m}(A_{v}A_{w},B_{v}B_{w}),\qquad
    Z_{w^{*\ell}}(m)=H_{m}(A_{w}^{\ell},B_{w}^{\ell}).
  \]
  Hence, if $v$, $w$, $v*w$ satisfy Hypothesis~\ref{hyp:positivity},
  and either $\lvert\theta_{2}\rvert\le1$, or $P, Q\in\ratQ$ and $d$ is
  not a square, then
  $\zeta_{v*w}(-m)=H_{m}(A_{v}A_{w},B_{v}B_{w})$, while
  $\zeta_{v}(-m)=H_{m}(A_{v},B_{v})$ and $\zeta_{w}(-m)=H_{m}(A_{w},B_{w})$.
\end{cor}
\begin{proof}
  The first formula is the definition of $Z_{w}(m;1)$ (take $n_{1}=1$)
  together with Lemma~\ref{lem:alpha}~(3). The others follow from
  Proposition~\ref{prop:lambda}~(1). The last sentence is
  Theorem~\ref{thm:negval}.
\end{proof}

\begin{thm}\label{thm:negrational}
  Assume $(C)$. Let $P, Q, w_{0}, w_{1}\in\ratQ$, and let $m\ge1$ with $\beta_{i}\neq1$
  for $0\le i\le m$ ($\beta_{i}$ as in~\eqref{eq:beta}). Then
  $Z_{w}(m)\in\ratQ$. Consequently, if
  $\lvert\theta_{2}\rvert\le1$ or $d$ is not a square, and $w$ satisfies
  Hypothesis~\ref{hyp:positivity}, then $\zeta_{w}(-m)\in\ratQ$.
\end{thm}
\begin{proof}
  If $d$ is a square in $\ratQ$, then $\sqrt{d}$, $\theta_{1}$,
  $\theta_{2}$, $A$, $B$ are rational, and so is $Z_{w}(m;1)$.

  Otherwise let $K=\ratQ(\sqrt{d})$, with nontrivial automorphism
  $\sigma$. Then $\sigma(\sqrt d)=-\sqrt d$, $\sigma(\theta_{1})=\theta_{2}$,
  $\sigma(A)=B$, $\sigma(B)=A$ (Proposition~\ref{prop:trace}~(1)), and
  $\sigma(\beta_{i'})=\beta_{m-i'}$.
  Write $Z_{w}(m;1)=\sum_{i'=0}^{m}T_{i'}$ with
  \[
    T_{i'}=(\sqrt{d})^{-m}\binom{m}{i'}(-1)^{i'}A^{m-i'}B^{i'}\frac{\beta_{i'}}{1-\beta_{i'}}\in K.
  \]
  Then
  \[
    \sigma(T_{i'})=(-1)^{m}(\sqrt{d})^{-m}\binom{m}{i'}(-1)^{i'}B^{m-i'}A^{i'}
    \frac{\beta_{m-i'}}{1-\beta_{m-i'}}=T_{m-i'},
  \]
  since $\binom{m}{i'}=\binom{m}{m-i'}$ and $(-1)^{m+i'}=(-1)^{m-i'}$.
  So $\sigma$ permutes the terms. Hence $\sigma(Z_{w}(m))=Z_{w}(m)$ and
  $Z_{w}(m)\in\ratQ$.
\end{proof}

\begin{remk}\label{remk:negval}
  \begin{enumerate}
  \item $A_{w}B_{w}=w_{1}^{2}-Pw_{0}w_{1}+Qw_{0}^{2}$, which equals
    $w_{0}^{2}f(w_{1}/w_{0})$ if $w_{0}\neq0$. If $d$ is not a square,
    this is $\mathrm{N}_{K/\ratQ}(A_{w})$, where $K=\ratQ(\sqrt d)$.
  \item The term $E_{w}(m)$ can be nonzero, and then $\zeta_{w}(-m)$ need
    not be rational. For $f(t)=t^{2}-2t-8$ (so $(P,Q)=(2,-8)$, which satisfies $(C)$) we have $\theta_{1}=4$,
    $\theta_{2}=-2$, and Kamano's Remark~2 notes that $s=-1$ is not a pole.
    For $w=\calF$ ($A=B=1$, $\lambda=1$) and $m=1$: $d=36$, $\beta_{0}=4$,
    $\beta_{1}=-2$, $\beta_{2}=4^{-1}(-2)^{2}=1$. Hence
    $Z_{\calF}(1)=\frac16\bigl(\frac{4}{1-4}-\frac{-2}{1+2}\bigr)=-\frac19$ and, with $i=2$,
    $E_{\calF}(1)=\frac16\cdot(-1)\cdot\frac{1!\,0!}{2!}\cdot\frac{1}{\log4}=-\frac{1}{12\log4}$,
    so $\zeta_{\calF}(-1)=-\frac19-\frac{1}{12\log4}$. This number is
    irrational: $\log4$ is irrational, since $e^{a}=4^{b}$ with integers
    $a, b\ge1$ would make $e$ algebraic, contradicting Hermite's theorem
    that $e$ is transcendental.
    The closed formula of Kamano~\cite[Theorem~7]{MR3064567} for
    $\zeta_{\calF}(-m)$ gives $-1/9$ here, so it misses the term
    $E_{\calF}(m)$ at such points. The example also shows that the
    rationality of the values at negative integers~\cite[Theorem~6.1]{MR4594520}
    does not extend to rational sequences whose minimal polynomial is
    separable but reducible, as suggested in~\cite[\S7]{MR4594520}. If
    $P, Q\in\ratQ$ and $f$ is irreducible over $\ratQ$, no such term occurs
    (Lemma~\ref{lem:alpha}~(4)).
  \item For $w=\calF$ with $f(t)=t^{2}-t-1$ we have $A=B=1$ and
    $\beta_{0}=\theta_{1}$, $\beta_{1}=\theta_{2}=-\theta_{1}^{-1}$, so
    \[
      \zeta_{\calF}(-1)=\frac{1}{\sqrt5}\Bigl(\frac{\theta_{1}}{1-\theta_{1}}
      -\frac{\theta_{2}}{1-\theta_{2}}\Bigr)
      =\frac{1}{\sqrt5}\Bigl(\frac{\theta_{1}}{\theta_{2}}-\frac{\theta_{2}}{\theta_{1}}\Bigr)
      =\frac{-\theta_{1}^{2}+\theta_{2}^{2}}{\sqrt5}=-1,
    \]
    using $1-\theta_{1}=\theta_{2}$, $1-\theta_{2}=\theta_{1}$,
    $\theta_{1}\theta_{2}=-1$, and
    $\theta_{1}^{2}-\theta_{2}^{2}=(\theta_{1}+\theta_{2})(\theta_{1}-\theta_{2})=\sqrt5$.
    Here $\lvert\theta_{2}\rvert<1$ and $\beta_{0},\beta_{1}\neq1$, so
    Theorem~\ref{thm:negval} applies.
  \end{enumerate}
\end{remk}

\subsection{Integrality of the values}
\label{sec:integrality}

Let $P, Q\in\intZ$ satisfy $(C)$, assume that $d$ is not a square, put
$K=\ratQ(\sqrt d)$, and let $m\ge1$ with $\beta_{i}\neq1$ for
$0\le i\le m$ ($\beta_{i}$ as in~\eqref{eq:beta}). Recall that
$\sqrt d=\theta_{1}-\theta_{2}>0$. The coefficients $c_{i}(m)$ of the
form $H_{m}$ of Corollary~\ref{cor:negform} contain $d^{-m/2}$. We
rewrite $H_{m}$ so that these powers of $\sqrt d$ cancel. This gives the
$p$-integrality of $Z_{w}(m)$ for every prime $p\nmid Q$ with
$p\nmid\mathrm{N}(1-\beta_{i})$, including the primes $p\mid d$. Put
$g_{i}=\beta_{i}/(1-\beta_{i})$ for $0\le i\le m$. For $k, n\ge0$ with
$k+n\le m$ let
\[
  \Delta^{n}g_{k}\coloneqq\sum_{l=0}^{n}\binom{n}{l}(-1)^{n-l}g_{k+l}
\]
be the $n$-th forward difference, and put
\begin{equation}
  \label{eq:gamma}
  \gamma_{j}\coloneqq(\sqrt d)^{-(m-j)}\,\Delta^{m-j}g_{j}\in K\qquad(0\le j\le m).
\end{equation}

\begin{lem}\label{lem:gamma}
  \begin{enumerate}
  \item As polynomials in $X$ and $Y$ over $K$,
    \[
      H_{m}(X,X-\sqrt d\,Y)=\sum_{j=0}^{m}(-1)^{m-j}\binom{m}{j}\gamma_{j}X^{m-j}Y^{j}.
    \]
    Hence, for $w\in\seqsp{f}{\realR}$,
    \[
      Z_{w}(m)=\sum_{j=0}^{m}(-1)^{m-j}\binom{m}{j}\gamma_{j}A_{w}^{m-j}w_{0}^{j},
    \]
    and $Z_{\calF}(m)=(-1)^{m}\gamma_{0}$.
  \item Let $p$ be a prime with $p\nmid Q$ and
    $p\nmid\mathrm{N}(1-\beta_{i})$ for $0\le i\le m$, and let
    $\mathfrak p$ be a prime ideal of $\mathcal O$ above $p$. Then
    $\gamma_{j}\in\mathcal O_{\mathfrak p}$ for $0\le j\le m$.
  \end{enumerate}
\end{lem}
\begin{proof}
  (1) By Corollary~\ref{cor:negform},
  $H_{m}(X,Y)=(\sqrt d)^{-m}\sum_{i=0}^{m}\binom{m}{i}(-1)^{i}g_{i}X^{m-i}Y^{i}$.
  Expand $(X-\sqrt d\,Y)^{i}=\sum_{j=0}^{i}\binom{i}{j}X^{i-j}(-\sqrt d\,Y)^{j}$.
  Since $\binom{m}{i}\binom{i}{j}=\binom{m}{j}\binom{m-j}{i-j}$, the
  coefficient of $X^{m-j}Y^{j}$ in $H_{m}(X,X-\sqrt d\,Y)$ is, with $i=j+l$,
  \[
    (\sqrt d)^{j-m}\binom{m}{j}\sum_{l=0}^{m-j}\binom{m-j}{l}(-1)^{l}g_{j+l}
    =(-1)^{m-j}\binom{m}{j}\gamma_{j}.
  \] By
  Proposition~\ref{prop:binet},
  $A_{w}-B_{w}=w_{0}(\theta_{1}-\theta_{2})=w_{0}\sqrt d$, so
  $B_{w}=A_{w}-\sqrt d\,w_{0}$, and $Z_{w}(m)=H_{m}(A_{w},B_{w})$ by
  Corollary~\ref{cor:negform}. For $w=\calF$ we have $w_{0}=0$ and
  $A_{\calF}=1$.

  (2) Since $\theta_{1}\theta_{2}=Q\notin\mathfrak p$, the elements
  $\theta_{1}$, $\theta_{2}$ of $\mathcal O$ are units of
  $\mathcal O_{\mathfrak p}$. If $1-\beta_{i}\in\mathfrak p$, then
  $\mathrm{N}(1-\beta_{i})=(1-\beta_{i})\sigma(1-\beta_{i})\in\mathfrak p\cap\intZ=p\intZ$;
  so $1-\beta_{i}$ is a unit of $\mathcal O_{\mathfrak p}$ for
  $0\le i\le m$. Fix $j$ and put $n=m-j$. If $n=0$, then
  $\gamma_{m}=g_{m}\in\mathcal O_{\mathfrak p}$. Let $n\ge1$. Put
  $\rho=\theta_{2}/\theta_{1}\in\mathcal O_{\mathfrak p}$,
  $\kappa=\theta_{1}^{m}$ and $x_{l}=\rho^{j+l}$ for $0\le l\le n$, so that
  $\beta_{j+l}=\kappa x_{l}$. By Lemma~\ref{lem:kamano}, $0<\lvert\rho\rvert<1$,
  so $x_{0},\dots,x_{n}$ are pairwise distinct. Let
  $\omega(x)=1/(1-\kappa x)$; then $g_{j+l}=\omega(x_{l})-1$. The divided
  differences of $\omega$ are defined by $\omega[x_{r}]=\omega(x_{r})$ and
  \[
    \omega[x_{r},\dots,x_{s}]=\frac{\omega[x_{r+1},\dots,x_{s}]-\omega[x_{r},\dots,x_{s-1}]}{x_{s}-x_{r}}
    \qquad(0\le r<s\le n).
  \]
  We show that
  \[
    \omega[x_{r},\dots,x_{s}]=\frac{\kappa^{s-r}}{\prod_{l=r}^{s}(1-\kappa x_{l})}
    \qquad(0\le r\le s\le n).
  \]
  This holds for $r=s$, and the induction step on $s-r$ is
  \begin{align*}
    \frac{\omega[x_{r+1},\dots,x_{s}]-\omega[x_{r},\dots,x_{s-1}]}{x_{s}-x_{r}}
    &=\frac{\kappa^{s-r-1}\bigl((1-\kappa x_{r})-(1-\kappa x_{s})\bigr)}{(x_{s}-x_{r})\prod_{l=r}^{s}(1-\kappa x_{l})}\\
    &=\frac{\kappa^{s-r}}{\prod_{l=r}^{s}(1-\kappa x_{l})}.
  \end{align*}
  These lie in $\mathcal O_{\mathfrak p}$, since
  $1-\kappa x_{l}=1-\beta_{j+l}$. Newton's interpolation formula gives
  \[
    \omega(x_{l})=\sum_{q=0}^{n}\omega[x_{0},\dots,x_{q}]\,\varpi_{q}(x_{l})
    \qquad(0\le l\le n),
  \]
  where $\varpi_{q}(x)=\prod_{r=0}^{q-1}(x-x_{r})$; the
  terms with $q>l$ vanish. The $n$-th difference of a constant is $0$, so
  \[
    \Delta^{n}g_{j}=\sum_{q=0}^{n}\omega[x_{0},\dots,x_{q}]
    \sum_{l=0}^{n}\binom{n}{l}(-1)^{n-l}\varpi_{q}(x_{l}).
  \]
  Write $\varpi_{q}(x)=\sum_{s}e_{s}x^{s}$ with $e_{s}\in\mathcal O_{\mathfrak p}$.
  Then
  $\sum_{l=0}^{n}\binom{n}{l}(-1)^{n-l}x_{l}^{s}=\rho^{js}(\rho^{s}-1)^{n}$,
  and $\rho^{s}-1=(\rho-1)(1+\rho+\dots+\rho^{s-1})\in(\rho-1)\mathcal O_{\mathfrak p}$
  (for $s=0$ it is $0$). Finally
  $\rho-1=(\theta_{2}-\theta_{1})/\theta_{1}=-\sqrt d/\theta_{1}$. Hence
  $\Delta^{n}g_{j}\in(\sqrt d)^{n}\mathcal O_{\mathfrak p}$, i.e.\
  $\gamma_{j}\in\mathcal O_{\mathfrak p}$.
\end{proof}

\begin{lem}\label{lem:coeffint}
  Let $p$ be a prime with $p\nmid d$ and $p\nmid\mathrm{N}(1-\beta_{i})$
  for $0\le i\le m$, and let $\mathfrak p$ be a prime ideal of
  $\mathcal O$ above $p$. Then the coefficients $c_{i}(m)$ of $H_{m}$
  lie in $\mathcal O_{\mathfrak p}$.
\end{lem}
\begin{proof}
  Since $P, Q\in\intZ$, the numbers $\theta_{1}$, $\theta_{2}$,
  $\sqrt d=\theta_{1}-\theta_{2}$ and $\beta_{i}$ lie in $\mathcal O$. The
  elements $\sqrt d$ and $1-\beta_{i}$ are units of $\mathcal O_{\mathfrak p}$,
  since $\mathrm{N}(\sqrt d)=-d$ and $p\nmid d$, and
  $p\nmid\mathrm{N}(1-\beta_{i})$. Hence
  $c_{i}(m)=(\sqrt d)^{-m}\binom{m}{i}(-1)^{i}\beta_{i}(1-\beta_{i})^{-1}\in\mathcal O_{\mathfrak p}$.
\end{proof}

Let $p$ be a prime. If $v\in\seqsp{f}{\ratQ}$ has
$v_{0}, v_{1}\in\intZ_{(p)}$, then $v_{n}\in\intZ_{(p)}$ for all
$n\ge0$, by the recurrence. For $h\ge0$, two such sequences $v$, $v'$
are \emph{congruent modulo $p^{h}$} if $v_{0}\equiv v'_{0}$ and
$v_{1}\equiv v'_{1}\pmod{p^{h}\intZ_{(p)}}$; then
$v_{n}\equiv v'_{n}$ for all $n\ge0$, again by the recurrence.

\begin{prop}\label{prop:residues}
  Let $p$ be a prime with $p\nmid\mathrm{N}(1-\beta_{i})$ for
  $0\le i\le m$, and assume that $p\nmid d$ or $p\nmid Q$. Let
  $v, v'\in\seqsp{f}{\ratQ}$ with $v_{0}, v_{1}, v'_{0}, v'_{1}\in\intZ_{(p)}$,
  and let $h\ge0$.
  \begin{enumerate}
  \item $Z_{v}(m)\in\intZ_{(p)}$.
  \item If $v$ and $v'$ are congruent modulo $p^{h}$, then
    $Z_{v}(m)\equiv Z_{v'}(m)\pmod{p^{h}}$.
  \item $Z_{\xi v}(m)=\xi^{m}Z_{v}(m)$ for $\xi\in\ratQ$. In particular, if
    $v_{0}\equiv0$ and $v_{1}\equiv\xi\pmod{p^{h}}$ with $\xi\in\intZ_{(p)}$,
    then $Z_{v}(m)\equiv\xi^{m}Z_{\calF}(m)\pmod{p^{h}}$.
  \end{enumerate}
\end{prop}
\begin{proof}
  Let $\mathfrak p$ be a prime ideal of $\mathcal O$ above $p$. The
  numbers $Z_{v}(m)$, $Z_{v'}(m)$ are rational by
  Theorem~\ref{thm:negrational}. Since
  $v_{0}, v_{1}\in\intZ_{(p)}\subset\mathcal O_{\mathfrak p}$, the numbers
  $A_{v}=v_{1}-v_{0}\theta_{2}$ and $B_{v}=v_{1}-v_{0}\theta_{1}$ lie in
  $\mathcal O_{\mathfrak p}$, and likewise for $v'$. Under the congruences
  of~(2), $A_{v}-A_{v'}$, $B_{v}-B_{v'}$ and $v_{0}-v'_{0}$ lie in
  $p^{h}\mathcal O_{\mathfrak p}$.

  If $p\nmid d$, then $Z_{v}(m)=H_{m}(A_{v},B_{v})$
  (Corollary~\ref{cor:negform}), and the coefficients of $H_{m}$ lie in
  $\mathcal O_{\mathfrak p}$ by Lemma~\ref{lem:coeffint}. If $p\nmid Q$,
  then Lemma~\ref{lem:gamma}~(1) gives
  $Z_{v}(m)=\sum_{j=0}^{m}(-1)^{m-j}\binom{m}{j}\gamma_{j}A_{v}^{m-j}v_{0}^{j}$
  with $\gamma_{j}\in\mathcal O_{\mathfrak p}$ by Lemma~\ref{lem:gamma}~(2).
  In both cases $Z_{v}(m)$ is a sum of terms in $\mathcal O_{\mathfrak p}$,
  and each term for $v$ is congruent to the one for $v'$ modulo
  $p^{h}\mathcal O_{\mathfrak p}$. Hence
  $Z_{v}(m)\in\mathcal O_{\mathfrak p}\cap\ratQ=\intZ_{(p)}$ and
  $Z_{v}(m)-Z_{v'}(m)\in p^{h}\mathcal O_{\mathfrak p}\cap\ratQ=p^{h}\intZ_{(p)}$.

  (3) We have $A_{\xi v}=\xi A_{v}$ and $B_{\xi v}=\xi B_{v}$, and $H_{m}$ is
  homogeneous of degree~$m$. For the last assertion apply~(2) to
  $v'=\xi\calF$, which has $v'_{0}=0$ and $v'_{1}=\xi$.
\end{proof}

\subsection{Reduction modulo \texorpdfstring{$p$}{p}}
\label{sec:quot-group-subgr}

We recall the reduction modulo $p$ of Laxton groups from
Aoki--Kida~\cite[\S5]{MR3902117}. Here $\bvec{w_{1}}{w_{0}}$ and $w^{*k}$
are as in Definition~\ref{def:laxton}, $\lambda_{w}=B_{w}/A_{w}$ is as
in~\eqref{eq:lambda}, and $\calF$ is the Lucas sequence
(Definition~\ref{def:lucasseq}). Let $P, Q\in\intZ$, $Q\neq0$, and let
$p$ be a prime with $p\nmid Q$.

Every class in $\GQ{f}$ has a representative $\pvec{w_{1}}{w_{0}}$ with
coprime $w_{0}, w_{1}\in\intZ$, unique up to sign. Let
$\projl{\finFq{p}}$ be the projective line over $\finFq{p}$, whose points are
the classes $\bvec{a}{b}$ of nonzero pairs $(a,b)\in\finFq{p}^{2}$ modulo
$\finFq{p}^{\times}$, and write $x\bmod p\in\finFq{p}$ for the residue of
$x\in\intZ$. The \emph{reduction map} is
\[
  \red{p}\colon\GQ{f}\to\projl{\finFq{p}},\qquad
  \bvec{w_{1}}{w_{0}}\mapsto\bvec{w_{1}\bmod p}{w_{0}\bmod p}.
\]
Let $\tilde{f}=f\bmod p$ (coefficientwise), and let
$\GFp{\tilde f}=\valsp{\tilde f}{\finFq{p}}^{\times}/\finFq{p}^{\times}\subset\projl{\finFq{p}}$.
Put
\[
  \Kfp{f}{p}\coloneqq\red{p}^{-1}(\GFp{\tilde f})\quad\text{and}\quad
  \Gfp{f}{p}\coloneqq\red{p}^{-1}\bigl(\bvec{1}{0}\bigr).
\]

\begin{thm}[Aoki--Kida {\cite[Theorem~9~(2)]{MR3902117}}]\label{thm:aokikida}
  Assume that $\tilde f$ is reducible over $\finFq{p}$ and $p\nmid d$, and
  let $\bar\theta_{1}, \bar\theta_{2}\in\finFq{p}$ be its roots. Then
  \[
    \Kfp{f}{p}/\Gfp{f}{p}\xrightarrow{\ \sim\ }\finFq{p}^{\times},\qquad
    \bvec{w_{1}}{w_{0}}\bmod\Gfp{f}{p}\mapsto
    (w_{1}-w_{0}\bar\theta_{1})(w_{1}-w_{0}\bar\theta_{2})^{-1},
  \]
  is an isomorphism.
\end{thm}

For odd $p\nmid d$, ``$\tilde f$ is reducible'' means $(d/p)=1$. For
$p\mid d$, Aoki and Kida show that $\Kfp{f}{p}/\Gfp{f}{p}$ is isomorphic
to the additive group $\finFq{p}$, by
$\bvec{w_{1}}{w_{0}}\mapsto-w_{0}(w_{1}-w_{0}\bar\theta)^{-1}$, where
$\bar\theta$ is the double root of $\tilde f$
\cite[Theorem~9~(2)]{MR3902117}. For $(d/p)=-1$, they show
$\Kfp{f}{p}/\Gfp{f}{p}\cong\mathbb F_{p^{2}}^{\times}/\finFq{p}^{\times}$,
a cyclic group of order $p+1$ \cite[Theorem~9~(1)]{MR3902117}.

In the rest of this section, $P, Q\in\intZ$, $d$ is not a square,
$K=\ratQ(\sqrt d)$, and $p$ is an odd prime; the assumption $p\nmid Q$
above is dropped unless we state it. Put $\delta_{p}=(d/p)\in\set{1,-1,0}$. For
$w\in\seqsp{f}{\ratQ}$ with $w_{0}, w_{1}\in\intZ$ put
$\eta_{w}=w_{1}-w_{0}P/2\in\intZ_{(p)}$ and
\begin{equation}
  \label{eq:tau}
  \tau_{w}\coloneqq
  \begin{cases}
    1&\text{if }\delta_{p}=1,\\
    \mathrm{N}(A_{w})&\text{if }\delta_{p}=-1,\\
    \eta_{w}&\text{if }\delta_{p}=0.
  \end{cases}
\end{equation}
In each case $p-\delta_{p}$ is the order of $\Kfp{f}{p}/\Gfp{f}{p}$
when $p\nmid Q$.

\begin{lem}\label{lem:liftp}
  Let $\Lambda$ be a commutative ring, $p$ a prime, and $l\ge1$. If $z\in \Lambda$
  satisfies $z\equiv1\pmod{p^{l}\Lambda}$, then $z^{p}\equiv1\pmod{p^{l+1}\Lambda}$.
  Consequently, $z\equiv1\pmod{p\Lambda}$ implies
  $z^{p^{e}}\equiv1\pmod{p^{e+1}\Lambda}$ for all $e\ge0$.
\end{lem}
\begin{proof}
  Write $z=1+p^{l}x$. Then
  $z^{p}=1+p^{l+1}x+\sum_{i=2}^{p-1}\binom{p}{i}p^{il}x^{i}+p^{pl}x^{p}$.
  For $2\le i\le p-1$, $p\mid\binom{p}{i}$, so the term lies in
  $p^{1+il}\Lambda\subset p^{l+1}\Lambda$. Also $pl\ge l+1$, because $p\ge2$ and $l\ge1$.
  The second assertion follows by induction on $e$.
\end{proof}

\begin{lem}\label{lem:laxtonpower}
  Let $w\in\seqsp{f}{\ratQ}$ with $w_{0}, w_{1}\in\intZ$, and assume $p\nmid\mathrm{N}(A_{w})$ if
  $\delta_{p}\neq0$. Let $e\ge0$, and put $u=w^{*(p-\delta_{p})p^{e}}$.
  Then $u_{0}, u_{1}\in\intZ$, and $u$ is congruent to
  $\tau_{w}^{p^{e}}\calF$ modulo $p^{e+1}$; that is,
  \[
    u_{0}\equiv0\quad\text{and}\quad u_{1}\equiv\tau_{w}^{p^{e}}\pmod{p^{e+1}}.
  \]
  If $\delta_{p}=0$, then
  $\mathrm{N}(A_{w})\equiv\eta_{w}^{2}\pmod{p\intZ_{(p)}}$; so
  $p\nmid\mathrm{N}(A_{w})$ implies that $\eta_{w}$ is a unit of
  $\intZ_{(p)}$.
\end{lem}
\begin{proof}
  Write $A=A_{w}$, $B=B_{w}$, $\nu=(p-\delta_{p})p^{e}$ and
  $\tau=\tau_{w}$, and let $\mathfrak p$ be a prime ideal of $\mathcal O$
  above $p$. We have $u_{0}, u_{1}\in\intZ$, since the formula of
  Lemma~\ref{lem:multiplication} has integer coefficients. By
  Proposition~\ref{prop:lambda}, $A_{u}=A^{\nu}$ and $B_{u}=B^{\nu}$; by
  Proposition~\ref{prop:binet}, $u_{0}\sqrt d=A_{u}-B_{u}$ and
  $u_{1}=A_{u}+u_{0}\theta_{2}$. Let $s_{0}\in\intZ_{(p)}$. We first show:
  \begin{equation}
    \label{eq:Apower}
    \text{if } A^{\nu}-s_{0},\ B^{\nu}-s_{0}\in p^{e+1}\sqrt d\,\mathcal O_{\mathfrak p},
    \text{ then } u_{0}\equiv0 \text{ and } u_{1}\equiv s_{0}\pmod{p^{e+1}}.
  \end{equation}
  Indeed, $u_{0}=\bigl((A^{\nu}-s_{0})-(B^{\nu}-s_{0})\bigr)/\sqrt d$ lies in
  $p^{e+1}\mathcal O_{\mathfrak p}\cap\ratQ=p^{e+1}\intZ_{(p)}$, and then
  $u_{1}-s_{0}=(A^{\nu}-s_{0})+u_{0}\theta_{2}\in p^{e+1}\mathcal O_{\mathfrak p}\cap\ratQ$.

  Let $\delta_{p}\neq0$. Since $\theta_{1}$ is a root of the monic
  polynomial $f\in\intZ[t]$ of discriminant $d$, we have
  $d=[\mathcal O:\intZ[\theta_{1}]]^{2}d_{K}$. As $p\nmid d$, $p$ divides
  neither $d_{K}$ nor the index, and $(d_{K}/p)=(d/p)=\delta_{p}$. In
  particular $\sqrt d$ is a unit of $\mathcal O_{\mathfrak p}$, so the
  modulus in~\eqref{eq:Apower} is $p^{e+1}\mathcal O_{\mathfrak p}$. By
  hypothesis $p\nmid\mathrm{N}(A)=\mathrm{N}(B)$
  (Proposition~\ref{prop:trace}~(1)), so $A$ and $B$ lie in no prime
  ideal above $p$.

  If $\delta_{p}=1$, then $\tilde f$ has two distinct roots, and $p$
  splits in $K$ by the Dedekind--Kummer theorem:
  $p\mathcal O=\mathfrak p\mathfrak p'$ with $\mathfrak p\neq\mathfrak p'$
  and $\mathcal O/\mathfrak p\cong\mathcal O/\mathfrak p'\cong\finFq{p}$.
  By Fermat's little theorem, $A^{p-1}\equiv1$ modulo $\mathfrak p$ and
  modulo $\mathfrak p'$. Since $\mathfrak p\cap\mathfrak p'=\mathfrak p\mathfrak p'=p\mathcal O$,
  we get $A^{p-1}\equiv1\pmod{p\mathcal O}$. Lemma~\ref{lem:liftp} with
  $z=A^{p-1}$ gives $A^{\nu}\equiv1\pmod{p^{e+1}\mathcal O}$, and likewise
  for $B$. Apply~\eqref{eq:Apower} with $s_{0}=1=\tau^{p^{e}}$.

  If $\delta_{p}=-1$, then $p$ is inert in $K$, since $p\nmid d_{K}$ and
  $(d_{K}/p)=-1$. So $\mathfrak p=p\mathcal O$, and the residue field
  $\mathcal O/p\mathcal O$ has $p^{2}$ elements. Since $p$ is unramified,
  the decomposition group $\mathrm{Gal}(K/\ratQ)$ maps isomorphically onto
  the Galois group of $\mathcal O/p\mathcal O$ over $\finFq{p}$. So
  $\sigma$ induces the Frobenius $x\mapsto x^{p}$, that is,
  $x^{p}\equiv\sigma(x)\pmod{p\mathcal O}$ for $x\in\mathcal O$. With
  $\sigma(A)=B$ this gives $A^{p+1}\equiv AB=\mathrm{N}(A)\pmod{p\mathcal O}$.
  Choose $\nu'\in\intZ$ with $\nu'\mathrm{N}(A)\equiv1\pmod{p^{e+1}}$. Then
  $z=\nu'A^{p+1}\equiv1\pmod{p\mathcal O}$, so
  $z^{p^{e}}\equiv1\pmod{p^{e+1}\mathcal O}$ by Lemma~\ref{lem:liftp}.
  Multiplying by $\mathrm{N}(A)^{p^{e}}$ and using
  $(\nu'\mathrm{N}(A))^{p^{e}}\equiv1\pmod{p^{e+1}}$ gives
  $A^{\nu}\equiv\mathrm{N}(A)^{p^{e}}\pmod{p^{e+1}\mathcal O}$. The same
  holds for $B$. Apply~\eqref{eq:Apower} with $s_{0}=\mathrm{N}(A)^{p^{e}}=\tau^{p^{e}}$.

  Let $\delta_{p}=0$, so $p\mid d$ and $\nu=p^{e+1}$. Since
  $\theta_{2}=(P-\sqrt d)/2$, we have $A=\eta_{w}+(w_{0}/2)\sqrt d$ and
  $B=\eta_{w}-(w_{0}/2)\sqrt d$. Hence
  $\mathrm{N}(A)=AB=\eta_{w}^{2}-w_{0}^{2}d/4\equiv\eta_{w}^{2}\pmod{p\intZ_{(p)}}$,
  which is the last assertion. As $p$ is odd, $w_{0}/2\in\mathcal O_{\mathfrak p}$,
  so $A\equiv\eta_{w}\equiv B\pmod{\sqrt d\,\mathcal O_{\mathfrak p}}$.
  Suppose $z\equiv b\pmod{p^{h}\sqrt d\,\mathcal O_{\mathfrak p}}$ with
  $z, b\in\mathcal O_{\mathfrak p}$ and $h\ge0$, and write
  $z=b+p^{h}\sqrt d\,y$. Then
  $z^{p}=b^{p}+\sum_{i=1}^{p}\binom{p}{i}b^{p-i}p^{hi}(\sqrt d)^{i}y^{i}$.
  For $1\le i\le p-1$, $p\mid\binom{p}{i}$, and the term lies in
  $p^{1+hi}\sqrt d\,\mathcal O_{\mathfrak p}\subset p^{h+1}\sqrt d\,\mathcal O_{\mathfrak p}$.
  For $i=p$, $(\sqrt d)^{p}=\sqrt d\,d^{(p-1)/2}\in p\sqrt d\,\mathcal O_{\mathfrak p}$,
  because $p\mid d$ and $p\ge3$; so that term lies in
  $p^{hp+1}\sqrt d\,\mathcal O_{\mathfrak p}\subset p^{h+1}\sqrt d\,\mathcal O_{\mathfrak p}$.
  Hence $z^{p}\equiv b^{p}\pmod{p^{h+1}\sqrt d\,\mathcal O_{\mathfrak p}}$.
  By induction on $h$, $A^{p^{h}}\equiv\eta_{w}^{p^{h}}\equiv B^{p^{h}}$
  modulo $p^{h}\sqrt d\,\mathcal O_{\mathfrak p}$ for every $h\ge0$.
  With $h=e+1$, \eqref{eq:Apower} gives $u_{0}\equiv0$ and
  $u_{1}\equiv\eta_{w}^{p^{e+1}}\pmod{p^{e+1}}$. Finally
  $\eta_{w}^{p^{e+1}}\equiv\eta_{w}^{p^{e}}\pmod{p^{e+1}\intZ_{(p)}}$: if
  $p\mid\eta_{w}$, both sides lie in $p^{p^{e}}\intZ_{(p)}$ and
  $p^{e}\ge e+1$; otherwise $\eta_{w}^{(p-1)p^{e}}\equiv1$, since
  $(\intZ/p^{e+1}\intZ)^{\times}$ has order $(p-1)p^{e}$.
\end{proof}

\begin{remk}\label{remk:AK}
  Assume $p\nmid Q$, and let $w\in\seqsp{f}{\ratQ}$ with coprime
  $w_{0}, w_{1}\in\intZ$ and $p\nmid\mathrm{N}(A_{w})$.
  \begin{enumerate}
  \item The class $\bvec{w_{1}}{w_{0}}$ lies in $\Kfp{f}{p}$, since
    $\mathrm{N}(A_{w})$ is the norm of $w_{1}-w_{0}t$ from $\ratQ[t]/(f)$,
    and its residue is the norm of $w_{1}-w_{0}t$ from
    $\finFq{p}[t]/(\tilde f)$.
  \item The case $e=0$ of Lemma~\ref{lem:laxtonpower} gives
    $u=w^{*(p-\delta_{p})}\equiv\tau_{w}\calF\pmod p$. Here
    $p\nmid\tau_{w}$: for $\delta_{p}=-1$ by assumption, and for
    $\delta_{p}=0$ by Lemma~\ref{lem:laxtonpower}. So $p\nmid u_{1}$, the
    greatest common divisor of $u_{0}$ and $u_{1}$ is prime to $p$, and the
    class of $u$ lies in $\Gfp{f}{p}$. This recovers the fact that the
    order of the class of $w$ in $\Kfp{f}{p}/\Gfp{f}{p}$ divides
    $p-\delta_{p}$, the order of this group. In addition, it gives the
    residue $\tau_{w}$ of $u_{1}$.
  \item Let $\delta_{p}=1$, let $\mathfrak p$ be the prime of $\mathcal O$
    above $p$ with $\theta_{i}\equiv\bar\theta_{i}\pmod{\mathfrak p}$, and
    identify $\mathcal O/\mathfrak p$ with $\finFq{p}$. Then
    $w_{1}-w_{0}\bar\theta_{1}\equiv B_{w}$ and
    $w_{1}-w_{0}\bar\theta_{2}\equiv A_{w}\pmod{\mathfrak p}$, and both are
    nonzero. So the isomorphism of Theorem~\ref{thm:aokikida} is
    $\lambda_{w}=B_{w}/A_{w}$ reduced modulo $\mathfrak p$.
  \end{enumerate}
\end{remk}

\subsection{Congruences in the Laxton exponent}
\label{sec:kummer-k}

Here $\beta_{i}$ is defined in~\eqref{eq:beta}, $Z_{w}(m)$ in
Lemma~\ref{lem:alpha}~(3), $H_{m}$ in Corollary~\ref{cor:negform},
$\delta_{p}$ and $\tau_{w}$ in~\eqref{eq:tau}, and the congruence of
sequences before Proposition~\ref{prop:residues}.

\begin{thm}\label{thm:congruence}
  Let $P, Q\in\intZ$ satisfy $(C)$, and assume that $d$ is not a square.
  Let $w_{0}, w_{1}\in\intZ$, and let $m\ge1$ with $\beta_{i}\neq1$ for
  $0\le i\le m$. Let $p$ be an odd prime such that
  \begin{enumerate}
  \item $p\nmid Q$ if $\delta_{p}=0$;
  \item $p\nmid \mathrm{N}(A_{w})=w_{1}^{2}-Pw_{0}w_{1}+Qw_{0}^{2}$;
  \item $p\nmid \mathrm{N}(1-\beta_{i})$ for $0\le i\le m$; if
    $\delta_{p}=0$, this is equivalent to $p\nmid 2^{m}-P^{m}$.
  \end{enumerate}
  Then:
  \begin{enumerate}
  \item[(i)] $Z_{w^{*\ell}}(m)\in\intZ_{(p)}$ for every $\ell\in\intZ$;
  \item[(ii)] if $e\ge0$ and $k\ge k'$ are integers with
    $k\equiv k'\pmod{(p-\delta_{p})p^{e}}$, and
    $\mu=(k-k')/((p-\delta_{p})p^{e})$, then
    \[
      Z_{w^{*k}}(m)\equiv\tau_{w}^{\mu mp^{e}}\,Z_{w^{*k'}}(m)\pmod{p^{e+1}}.
    \]
  \end{enumerate}
  If moreover $w^{*k}$ and $w^{*k'}$ satisfy
  Hypothesis~\ref{hyp:positivity}, the same congruence holds with
  $Z_{w^{*k}}(m)$, $Z_{w^{*k'}}(m)$ replaced by $\zeta_{w^{*k}}(-m)$,
  $\zeta_{w^{*k'}}(-m)$. If $\delta_{p}=0$, then (i) and (ii) hold without~(2),
  with $k, k', \ell\ge0$.
\end{thm}
\begin{proof}
  Write $A=A_{w}$.

  We prove the equivalence in~(3). Let $\delta_{p}=0$, and let
  $\mathfrak p$ be a prime ideal of $\mathcal O$ above $p$. Since
  $d\in p\intZ\subset\mathfrak p$ and $(\sqrt d)^{2}=d$, we have
  $\sqrt d\in\mathfrak p$. Hence
  $2\theta_{1}=P+\sqrt d\equiv P$ and $2\theta_{2}\equiv P\pmod{\mathfrak p}$,
  so $2^{m}(1-\beta_{i})\equiv2^{m}-P^{m}\pmod{\mathfrak p}$ for
  $0\le i\le m$, and $2\notin\mathfrak p$. If $p\mid2^{m}-P^{m}$, then
  $1-\beta_{0}\in\mathfrak p$, so
  $\mathrm{N}(1-\beta_{0})\in\mathfrak p\cap\intZ=p\intZ$. If
  $p\nmid2^{m}-P^{m}$, then $1-\beta_{i}\notin\mathfrak p$ for every
  $\mathfrak p$ above $p$; as $\sigma(1-\beta_{i})=1-\beta_{m-i}$, also
  $\mathrm{N}(1-\beta_{i})\notin\mathfrak p$, i.e.\
  $p\nmid\mathrm{N}(1-\beta_{i})$.

  If $\delta_{p}\neq0$, then $p\nmid d$; if $\delta_{p}=0$, then
  $p\nmid Q$ by~(1). With~(3), Proposition~\ref{prop:residues} applies.

  (i) By Lemma~\ref{lem:multiplication},
  $w*\pvec{w_{1}-Pw_{0}}{-w_{0}}=\mathrm{N}(A)\calF$. By~(2),
  $\mathrm{N}(A)$ is a unit of $\intZ_{(p)}$, so $w$ is invertible and
  $w^{*(-1)}=\mathrm{N}(A)^{-1}\pvec{w_{1}-Pw_{0}}{-w_{0}}$ has terms
  $0$ and $1$ in $\intZ_{(p)}$. By the formula of
  Lemma~\ref{lem:multiplication}, if $v_{0}, v_{1}, v'_{0}, v'_{1}\in\intZ_{(p)}$,
  then the terms $0$ and $1$ of $v*v'$ lie in $\intZ_{(p)}$. Hence
  $(w^{*\ell})_{0}, (w^{*\ell})_{1}\in\intZ_{(p)}$ for every
  $\ell\in\intZ$. Proposition~\ref{prop:residues}~(1) gives~(i).

  (ii) Put $u=w^{*(p-\delta_{p})p^{e}}$ and $\xi=\tau_{w}^{p^{e}}$. By
  Lemma~\ref{lem:laxtonpower}, $u$ is congruent to $\xi\calF$ modulo
  $p^{e+1}$. The terms $0$ and $1$ of $v*u$ are polynomials with integer
  coefficients in $v_{0}$, $v_{1}$, $u_{0}$, $u_{1}$
  (Lemma~\ref{lem:multiplication}). So, for
  $v\in\seqsp{f}{\ratQ}$ with $v_{0}, v_{1}\in\intZ_{(p)}$, $v*u$ is
  congruent to $v*(\xi\calF)=\xi v$ modulo $p^{e+1}$. Applying this $\mu$
  times, starting with $v=w^{*k'}$, shows that $w^{*k}=w^{*k'}*u^{*\mu}$
  is congruent to $\xi^{\mu}w^{*k'}$ modulo $p^{e+1}$. By
  Proposition~\ref{prop:residues}~(2) and~(3),
  \[
    Z_{w^{*k}}(m)\equiv Z_{\xi^{\mu}w^{*k'}}(m)=\xi^{\mu m}Z_{w^{*k'}}(m)\pmod{p^{e+1}}.
  \]

  If $\delta_{p}=0$ and $k, k', \ell\ge0$, then the terms $0$ and $1$ of
  $w^{*\ell}$ are integers, and Lemma~\ref{lem:laxtonpower} does not need~(2); so the argument
  above holds without~(2). Since $d$ is not a square,
  Theorem~\ref{thm:negval} gives $\zeta_{w^{*\ell}}(-m)=Z_{w^{*\ell}}(m)$
  whenever $w^{*\ell}$ satisfies Hypothesis~\ref{hyp:positivity}. This
  proves the last assertions.
\end{proof}

For $\delta_{p}=1$ the factor $\tau_{w}^{\mu mp^{e}}$ is $1$, and for
$\delta_{p}=-1$ it is $\mathrm{N}(A_{w})^{\mu mp^{e}}$. For
$\delta_{p}=0$ we have $k-k'=\mu p^{e+1}$, and the factor may be
replaced by $\eta_{w}^{(k-k')m}$, since
$\eta_{w}^{p^{e+1}}\equiv\eta_{w}^{p^{e}}\pmod{p^{e+1}\intZ_{(p)}}$ (end of
the proof of Lemma~\ref{lem:laxtonpower}).

\begin{cor}\label{cor:fermat}
  In the situation of Theorem~\ref{thm:congruence}, assume that
  $u=w^{*(p-\delta_{p})}$ satisfies Hypothesis~\ref{hyp:positivity}. Then
  $\zeta_{u}(-m)$ and $\zeta_{\calF}(-m)$ lie in $\intZ_{(p)}$, and
  \[
    \zeta_{w^{*(p-\delta_{p})}}(-m)\equiv\tau_{w}^{m}\,\zeta_{\calF}(-m)\pmod{p}.
  \]
\end{cor}
\begin{proof}
  Take $e=0$, $k=p-\delta_{p}$ and $k'=0$ in
  Theorem~\ref{thm:congruence}; then $\mu=1$ and $w^{*0}=\calF$. The
  sequence $\calF$ satisfies Hypothesis~\ref{hyp:positivity}
  (Lemma~\ref{lem:positivity}, $A_{\calF}=1$).
\end{proof}

Theorem~\ref{thm:congruence} assumes $p\nmid Q$ only if $\delta_{p}=0$.
If $p\mid Q$, then $d\equiv P^{2}\pmod p$, so $\delta_{p}\in\set{1,0}$;
the comparison with Aoki--Kida in Remark~\ref{remk:AK} needs $p\nmid Q$.

\begin{remk}\label{remk:fermathyp}
  Hypothesis~(3) of Theorem~\ref{thm:congruence} cannot be omitted, even
  in Corollary~\ref{cor:fermat}.
  \begin{enumerate}
  \item For $f(t)=t^{2}-t-1$, $p=11$ and $m=5$, we have
    $\mathrm{N}(1-\beta_{0})=\mathrm{N}(1-\theta_{1}^{5})=-11$, and
    $\zeta_{\calF}(-5)=-7/22$ is not $11$-integral.
  \item For $f(t)=t^{2}-3t-1$ we have $d=13$. Take $p=3$, so
    $\delta_{p}=(13/3)=1$, and $m=2$. Then
    $\mathrm{N}(1-\beta_{0})=-9$, so (3) fails. Take
    $w=\pvec{3}{1}$; then $\mathrm{N}(A_{w})=-1$, so (2) holds. Since
    $A_{w}>0$ and $-1<\lambda_{w}<0$, the sequence
    $w^{*2}=\pvec{10}{3}$, with terms $10, 33, 109,\dots$, satisfies
    Hypothesis~\ref{hyp:positivity} by
    Proposition~\ref{prop:powerpositive}. Both values are $3$-integral,
    but $\zeta_{w^{*2}}(-2)=-10\not\equiv0=\zeta_{\calF}(-2)\pmod3$.
  \end{enumerate}
\end{remk}

\begin{cor}\label{cor:period}
  Assume the hypotheses of Theorem~\ref{thm:congruence}, including~(2)
  when $\delta_{p}=0$. Put $M=p-1$ if
  $\delta_{p}=1$, $M=p^{2}-1$ if $\delta_{p}=-1$, and $M=p(p-1)$ if
  $\delta_{p}=0$. If $e\ge0$ and $k\equiv k'\pmod{Mp^{e}}$, then
  \[
    Z_{w^{*k}}(m)\equiv Z_{w^{*k'}}(m)\pmod{p^{e+1}}.
  \]
\end{cor}
\begin{proof}
  We may assume $k\ge k'$. If $\delta_{p}=1$, this is
  Theorem~\ref{thm:congruence}, whose factor is then $1$. Otherwise
  $M=(p-\delta_{p})(p-1)$, so $\mu=(p-1)\mu_{1}$ for an integer $\mu_{1}$,
  and the factor is $\bigl(\tau_{w}^{(p-1)p^{e}}\bigr)^{m\mu_{1}}$. Here
  $\tau_{w}$ is a unit of $\intZ_{(p)}$, by~(2) and, if $\delta_{p}=0$,
  Lemma~\ref{lem:laxtonpower}. Since $(\intZ/p^{e+1}\intZ)^{\times}$ has
  order $(p-1)p^{e}$, we have $\tau_{w}^{(p-1)p^{e}}\equiv1\pmod{p^{e+1}}$.
\end{proof}

\begin{cor}\label{cor:padicinterp}
  In the situation of Corollary~\ref{cor:period}, fix
  $a\in\set{0,\dots,M-1}$ and put $\mathcal I_{a}(n)=Z_{w^{*(a+Mn)}}(m)$ for
  $n\in\intZ$. Then $\mathcal I_{a}$ extends uniquely to a continuous function
  $\intZ_{p}\to\intZ_{p}$. It satisfies
  $\lvert \mathcal I_{a}(x)-\mathcal I_{a}(y)\rvert_{p}\le p^{-1}\lvert x-y\rvert_{p}$.
\end{cor}
\begin{proof}
  Let $n, n'\in\intZ$ with $\lvert n-n'\rvert_{p}=p^{-e}$, i.e.\
  $p^{e}\parallel n-n'$. Then $Mp^{e}\mid M(n-n')$, so by
  Corollary~\ref{cor:period},
  $\lvert \mathcal I_{a}(n)-\mathcal I_{a}(n')\rvert_{p}\le p^{-(e+1)}=p^{-1}\lvert n-n'\rvert_{p}$.
  Thus $\mathcal I_{a}$ is Lipschitz, hence uniformly continuous, on the dense
  subset $\intZ$ of $\intZ_{p}$. Its values lie in the complete space
  $\intZ_{p}$, so it extends uniquely to a continuous map on $\intZ_{p}$, and
  the Lipschitz bound persists by continuity.
\end{proof}

\begin{remk}\label{remk:kummerk}
  The proof of Theorem~\ref{thm:congruence} uses only
  Lemma~\ref{lem:laxtonpower}, the bilinearity of the Laxton product, and
  Proposition~\ref{prop:residues} for the fixed $m$. Accordingly,
  Corollary~\ref{cor:padicinterp} is an interpolation in the Laxton
  exponent $k$, not in $m$. Let $p\nmid Q$, and let $w_{0}, w_{1}\in\intZ$
  with $p\nmid\mathrm{N}(A_{w})$. Put $u=w^{*(p-\delta_{p})}$, and let
  $g$ be the greatest common divisor of $u_{0}$ and $u_{1}$. Then
  $g^{2}$ divides $\mathrm{N}(A_{u})=\mathrm{N}(A_{w})^{p-\delta_{p}}$, so
  $p\nmid g$. The description of Aoki and Kida alone shows that the class
  of $u$ lies in $\Gfp{f}{p}$, i.e.\ $u\equiv\xi\calF\pmod p$ for some
  $\xi\in\intZ$ prime to $p$. With Proposition~\ref{prop:residues}~(3),
  this gives Corollary~\ref{cor:fermat} for the values $Z$, with the
  factor $\xi^{m}$ in place of $\tau_{w}^{m}$.
  Lemma~\ref{lem:laxtonpower} identifies $\xi\equiv\tau_{w}\pmod p$ and
  lifts the congruence to the modulus $p^{e+1}$. Theorem~\ref{thm:congruence}
  also gives the $p$-integrality~(i) for $p\mid d$, which rests on
  Lemma~\ref{lem:gamma}~(2).
\end{remk}

\begin{remk}\label{remk:interppos}
  Corollary~\ref{cor:padicinterp} concerns the values $Z_{w^{*k}}(m)$
  for all $k\in\intZ$, and needs no positivity. In particular, the
  continuous extension of the values for $k\ge0$ takes the values
  $Z_{w^{*k}}(m)$ at the negative integers $k$. By
  Theorem~\ref{thm:negval}, $Z_{w^{*k}}(m)=\zeta_{w^{*k}}(-m)$ when
  $w^{*k}$ satisfies Hypothesis~\ref{hyp:positivity}. If
  $A_{w}>0$, then by Proposition~\ref{prop:powerpositive} this holds for
  all $k\ge0$ if $\lvert\lambda_{w}\rvert\le1$, and fails for all large
  even $k$ if $\lvert\lambda_{w}\rvert>1$. If $A_{w}<0$, it fails for every
  odd $k$, since $A_{w^{*k}}=A_{w}^{k}<0$. For $k<0$, apply this to
  $w^{*(-1)}$, for which $A_{w^{*(-1)}}=A_{w}^{-1}$ and
  $\lambda_{w^{*(-1)}}=\lambda_{w}^{-1}$: if $A_{w}>0$, the hypothesis
  holds for all $k<0$ if $\lvert\lambda_{w}\rvert\ge1$, and fails for all
  even $k<0$ with large $\lvert k\rvert$ if $\lvert\lambda_{w}\rvert<1$.
  In particular, the sequence $u=w^{*(p-\delta_{p})}$ of
  Corollary~\ref{cor:fermat} satisfies Hypothesis~\ref{hyp:positivity} if
  $A_{w}>0$ and $\lvert\lambda_{w}\rvert\le1$.
\end{remk}

\begin{ex}\label{ex:ramified}
  We computed the values below exactly from Corollary~\ref{cor:negform}
  with PARI/GP. Let $f=t^{2}-t-1$, so $d=5$, and $p=5$, so $\delta_{p}=0$.
  Take $w=\pvec{5}{2}$, i.e.\
  $w_{0}=2$, $w_{1}=5$; then $A_{w}=4+\sqrt5$, $B_{w}=4-\sqrt5$,
  $\mathrm{N}(A_{w})=11$, $\tau_{w}=\eta_{w}=4$, and
  $0<\lambda_{w}=(4-\sqrt5)/(4+\sqrt5)<1$. By
  Proposition~\ref{prop:powerpositive}, every $w^{*k}$ with $k\ge0$
  satisfies Hypothesis~\ref{hyp:positivity}. Here $w^{*5}=\pvec{6829}{4210}$, and
  $4210\equiv0$, $6829\equiv4\pmod5$, as in
  Lemma~\ref{lem:laxtonpower}. For $m=1$ and $m=3$ the hypotheses of
  Theorem~\ref{thm:congruence} hold ($\mathrm{N}(A_{w})=11$, $2-1=1$, $2^{3}-1=7$, and $\beta_{i}\neq1$
  by Lemma~\ref{lem:alpha}~(4)). We have $\zeta_{\calF}(-1)=-1$ and
  $\zeta_{w^{*5}}(-1)=-11039\equiv4\cdot(-1)\pmod5$, while
  $-11039\not\equiv\zeta_{\calF}(-1)$ and
  $-11039\not\equiv\mathrm{N}(A_{w})\zeta_{\calF}(-1)\pmod5$. Also
  $\zeta_{\calF}(-3)=1/2$ and $\zeta_{w^{*5}}(-3)=-195912572041/2\equiv4^{3}\cdot\frac12\pmod5$.
\end{ex}

\section{The case of units of real quadratic fields}
\label{sec:units}

In this section $P\in\intZ$ and $Q\in\set{\pm1}$, with $P\ge1$ if $Q=-1$
and $P\ge3$ if $Q=1$. Recall $d=P^{2}-4Q$, the roots $\theta_{1}$,
$\theta_{2}$ and the numbers $A_{w}$, $B_{w}$
(Proposition~\ref{prop:binet}), the Laxton group and its product $*$
(Definition~\ref{def:laxton}), $\lambda_{w}$~\eqref{eq:lambda}, and the
condition $(C)$ (Definition~\ref{def:condc}). Then $d$ is not a square
(Proposition~\ref{prop:units}~(1) below). Let $K=\ratQ(\sqrt d)$, with nontrivial automorphism
$\sigma$, ring of integers $\mathcal O$, and fundamental unit
$\varepsilon>1$ (the generator $>1$ of $\mathcal O^{\times}/\set{\pm1}$).
Put $K^{1}=\set{x\in K^{\times} | x\sigma(x)=1}$.

\begin{prop}\label{prop:units}
  \begin{enumerate}
  \item The condition $(C)$ holds, $d$ is not a square, and
    $\theta_{1}=\varepsilon^{g}$ for some integer $g\ge1$. Moreover
    $\mathrm{N}(\varepsilon)^{g}=Q$.
  \item The map $\bvec{w_{1}}{w_{0}}\mapsto A_{w}\ratQ^{\times}$ is a group
    isomorphism $\GQ{f}\to K^{\times}/\ratQ^{\times}$. For every
    $\alpha\in K^{\times}$ there is $w\in\seqsp{f}{\ratQ}$ with
    \[
      w_{n}=\mathrm{Tr}_{K/\ratQ}(\alpha\,\varepsilon^{gn})\qquad(n\in\intZ),
    \]
    namely the one with $A_{w}=\alpha\sqrt d$; and every $w$ with
    $A_{w}\neq0$ arises in this way, with $\alpha=A_{w}/\sqrt d$.
  \item The map $\bvec{w_{1}}{w_{0}}\mapsto\lambda_{w}=\sigma(A_{w})/A_{w}$~\eqref{eq:lambda}
    is a group isomorphism $\GQ{f}\to K^{1}$.
  \end{enumerate}
  The first assertion of~(2) is the case $R=\ratQ$ of the description of
  $G_{R}(f)$ for irreducible $f$ by Aoki and Kida~\cite{MR3902117}.
\end{prop}
\begin{proof}
  (1) If $Q=-1$ and $P\ge1$, then $Q<P-1$; if $Q=1$ and $P\ge3$, then
  $Q\le P-1$ and $P>2$. So $(C)$ holds. If $Q=-1$, then $d=P^{2}+4$ lies
  strictly between $P^{2}$ and $(P+2)^{2}$ and has the parity of $P^{2}$, so
  it is not a square. If $Q=1$, then $(P-1)^{2}<P^{2}-4<P^{2}$ for $P\ge3$.
  Hence $f$ is irreducible over $\ratQ$ and $\theta_{1}\in K$. It is a root
  of the monic $f\in\intZ[t]$, so $\theta_{1}\in\mathcal O$, and
  $\theta_{1}\sigma(\theta_{1})=\theta_{1}\theta_{2}=Q=\pm1$, so
  $\theta_{1}\in\mathcal O^{\times}$. By Lemma~\ref{lem:kamano},
  $\theta_{1}>1$. Since $\mathcal O^{\times}=\set{\pm\varepsilon^{i} | i\in\intZ}$
  and $\varepsilon>1$, the only units $>1$ are $\varepsilon^{i}$ with
  $i\ge1$. Taking norms gives $\mathrm{N}(\varepsilon)^{g}=Q$.

  (2) Since $f$ is irreducible, $\ratQ[t]/(f(t))\to K$, $t\mapsto\theta_{2}$,
  is a field isomorphism; under \eqref{eq:three-modules} it sends $w$ to
  $A_{w}$. Hence $\valsp{f}{\ratQ}^{\times}\cong K^{\times}$, and dividing by
  $\ratQ^{\times}$ gives the first assertion. The formula for $w_{n}$ is
  Proposition~\ref{prop:trace}~(1) with $\theta_{1}=\varepsilon^{g}$. Given
  $\alpha\in K^{\times}$, write $\alpha\sqrt d=a-b\theta_{2}$ with
  $a, b\in\ratQ$ (possible since $1, \theta_{2}$ is a basis of $K$ over
  $\ratQ$) and take $(w_{1},w_{0})=(a,b)$.

  (3) By (2) and since $\sigma(A_{w})=B_{w}$
  (Proposition~\ref{prop:trace}~(1)), the map is the composite of
  $\GQ{f}\cong K^{\times}/\ratQ^{\times}$ with
  $K^{\times}/\ratQ^{\times}\to K^{1}$, $x\mapsto\sigma(x)/x$. The latter is
  well defined and a homomorphism. It is injective, since $\sigma(x)=x$
  means $x\in\ratQ^{\times}$. It is surjective by Hilbert's Theorem~90 for
  the cyclic extension $K/\ratQ$.
\end{proof}

\begin{remk}\label{remk:units}
  \begin{enumerate}
  \item Proposition~\ref{prop:units} identifies the sequences
    $\mathrm{Tr}(\alpha\varepsilon^{gn})$, $\alpha\in K^{\times}$, with the
    elements of $\GQ{f}$. The Laxton product corresponds to the product of
    the elements $\alpha\sqrt d$, and $\lambda$ to $\alpha\mapsto\sigma(\alpha)/\alpha$
    up to the sign $\sigma(\sqrt d)/\sqrt d=-1$. In these terms,
    Corollary~\ref{cor:zetaprod} and Corollary~\ref{cor:negform} express
    $\zeta_{v*w}$ and $\zeta_{v*w}(-m)$ through the products of the
    corresponding elements of $K^{\times}$.
  \item The sign of $Q=\mathrm{N}(\varepsilon)^{g}$ governs positivity. If $Q=-1$,
    $v_{0},w_{0}\ge0$ and $v_{1},w_{1}>0$, then the hypothesis on $v*w$ in
    Corollary~\ref{cor:zetaprod} follows from
    Proposition~\ref{prop:prodpositive}. If $Q=1$, it does not; see
    Remark~\ref{remk:counterexamples}~(2), where $f(t)=t^{2}-3t+1$ and
    $\theta_{1}=\varepsilon^{2}$ for $K=\ratQ(\sqrt5)$.
  \item By Lemma~\ref{lem:alpha}~(4), $s=-m$ is a pole of $\zeta_{w}$
    exactly when $Q^{m/2}=1$, that is, for even $m$ if $Q=1$ and for
    $m\equiv0\pmod4$ if $Q=-1$; the residue is nonzero since
    $\lambda_{w}\neq0$. For $w=\calF$ this is
    Kamano~\cite[Proposition~6]{MR3064567}; for instance, $\zeta_{\calF}$
    has a pole at $s=-4$ for $f(t)=t^{2}-t-1$. So in the unit case Section~\ref{sec:special-values-lucas}
    covers only the remaining $m$.
  \item In Theorem~\ref{thm:congruence},
    the condition $p\nmid Q$ of~(1) holds automatically. Primes $p\mid d$ occur,
    for instance $p=5$ for $f(t)=t^{2}-t-1$
    (Example~\ref{ex:ramified}). With $A_{w}=\alpha\sqrt d$ and
    $\delta_{p}\neq0$, the condition~(2), $p\nmid \mathrm{N}(A_{w})$, reads
    $p\nmid \mathrm{N}_{K/\ratQ}(\alpha)$ when $\alpha\in\mathcal O$, since
    then $p\nmid d$.
  \end{enumerate}
\end{remk}

\section{Questions not treated here}
\label{sec:conclusion}

We list some questions that this paper does not treat.
\begin{enumerate}
\item The values of $\zeta_{w}$ at the positive integers: representations
  by Lambert series and theta functions in the manner of Landau, and their
  irrationality and transcendence (cf.~Murty~\cite{MR3156859}).
\item Congruences in $m$ of Kummer type and a $p$-adic interpolation of
  $\zeta_{w}(-m)$ in the variable $m$. We will treat these in a separate
  paper.
\item The trivial zeros of $\zeta_{w}$ for general $w$ (for $w=\calF$
  see~\cite[Corollaries~8 and~12]{MR3064567}).
\item The prime $p=2$ in the congruences of
  \S\ref{sec:special-values-lucas}. Proposition~\ref{prop:residues} includes
  $p=2$, but Theorem~\ref{thm:congruence} assumes $p$ odd. If $2\mid d$ and $2\nmid Q$, then
  $P$ is even, $\sqrt d\in2\mathcal O$, and $\theta_{1}\equiv\theta_{2}$,
  $\theta_{1}^{2}\equiv Q\equiv1$ modulo every prime ideal $\mathfrak p$
  above $2$; so $\theta_{1}\equiv\theta_{2}\equiv1$ and
  $2\mid\mathrm{N}(1-\beta_{i})$ for all $i$. In this case $Z_{w}(m)$ need
  not be $2$-integral: $Z_{\calF}(1)=-1/2$ for $f(t)=t^{2}-2t-1$.
\item Analogs for other unit groups of rank one: cubic fields that are
  not totally real (for a family of such fields with explicit
  fundamental units, see~\cite{KimuraUmemoto2026}), totally imaginary quartic fields, and the elements of an
  imaginary quadratic field that are units outside one fixed prime ideal.
\end{enumerate}

\section*{Acknowledgments}

The author thanks M.~Kida for drawing his attention to the work of
Aoki and Kida~\cite{MR3902117}, and G.~Shibukawa for pointing out the relation between the
proof of Theorem~\ref{thm:expforseq} and the work of Kaneko, Kurokawa and
Wakayama~\cite{MR2069738}.

\section*{Declarations}

Declarations of interest: none.
This research did not receive any specific grant from funding agencies
in the public, commercial, or not-for-profit sectors.

\section*{\AIdeclarationtitle}

The author used Claude (Anthropic) to write and run PARI/GP scripts for
numerical and exact checks of the statements, to search the literature,
and to improve the English. The mathematical content and the proofs are
the author's. The author checked all output and takes full responsibility
for the content of this paper.

\bibliography{paper1}

\end{document}